\documentclass[brochure,12pt]{bourbaki}
\usepackage[matrix,arrow]{xy}
\usepackage{amssymb,amsfonts,amsmath,footnote}
\usepackage{graphicx,color,epsfig} 
\usepackage[francais]{babel}
\newcommand{\C}{\mathbb{C}}
\newcommand{\N}{\mathbb{N}}
\renewcommand{\P}{\mathbb{P}}
\newcommand{\R}{\mathbb{R}}
\newcommand{\Z}{\mathbb{Z}}

\newcommand{\cB}{{\mathcal{B}}}
\newcommand{\cE}{{\mathcal{E}}}
\newcommand{\cJ}{{\mathcal{J}}}
\newcommand{\cM}{{\mathcal{M}}}
\newcommand{\cN}{{\mathcal{N}}}
\newcommand{\cO}{{\mathcal{O}}}
\newcommand{\cP}{{\mathcal{P}}}
\newcommand{\cR}{{\mathcal{R}}}

\newcommand{\cW}{{\mathcal{W}}}

\newcommand{\eps}{\varepsilon}

\newcommand{\oD}{{\overline{D}}}
\newcommand{\ux}{\underline{x}}
\newcommand{\uz}{{\underline{z}}}

\newcommand{\pp}{{\mathfrak{p}}}
\renewcommand{\ss}{{\mathfrak{s}}}

\newcommand{\p}{\mathbf{p}}
\newcommand{\q}{\mathbf{q}}

\newcommand{\dbar}{{\bar\partial}}

\newcommand{\coker}{{\mathrm{coker}\,}}
\newcommand{\ev}{{\mathrm{ev}}}
\newcommand{\oev}{{\overline{\mathrm{ev}}}}
\newcommand{\ind}{{\mathrm{ind}}}

\newcommand{\Pin}{{\mathrm{Pin}}}
\newcommand{\Spin}{{\mathrm{Spin}}}

\date{Avril 2011}
\bbkannee{63\`eme ann\'ee, 2010-2011}
\bbknumero{1036}
\title{Invariants de Welschinger}
\subtitle{}
\author{Alexandru OANCEA}
\address{Institut de Recherche Math\'ematique Avanc\'ee (IRMA)\\
CNRS et Universit\'e de Strasbourg\\
7, Rue Ren\'e Descartes\\
F--67084 Strasbourg Cedex}
\email{oancea@math.unistra.fr}

\begin{document}
\maketitle



Le but de cet expos\'e est de pr\'esenter des invariants d\'ecouverts par Welschinger qui sont adapt\'es \`a des probl\`emes de g\'eom\'etrie \'enum\'erative r\'eelle.  Ces probl\`emes \'enum\'eratifs sont classiquement formul\'es dans le cadre de la g\'eom\'etrie alg\'ebrique r\'eelle mais ils trouvent leur solution la plus naturelle dans le cadre de la g\'eom\'etrie symplectique r\'eelle. Ceci permet en particulier de les \'etudier via des techniques sp\'ecifiques puissantes comme la th\'eorie symplectique des champs. Les invariants de Welschinger sont des analogues r\'eels de certains invariants de Gromov-Witten.

\section{Introduction}

\bigskip

Une \emph{vari\'et\'e symplectique r\'eelle} $(X,\omega,c_X)$ est une vari\'et\'e diff\'erentiable munie d'une $2$-forme ferm\'ee non-d\'eg\'en\'er\'ee $\omega$ -- la forme symplectique -- et d'une involution anti-symplectique $c_X:X\to X$, v\'erifiant $c_X^*\omega=$\,\,-$\omega$ -- la \emph{structure r\'eelle}. La dimension de $X$ est n\'ecessairement paire, not\'ee $2n$. S'il est non-vide, le \emph{lieu r\'eel} $\R X=\mathrm{Fix}(c_X)\subset X$ est une sous-vari\'et\'e lisse lagrangienne, i.e. $\dim\, \R X=n$ et $\omega|_{\R X}=0$. Ceci d\'ecoule de l'\'enonc\'e analogue pour les structures r\'eelles lin\'eaires sur $\R^{2n}$ et du th\'eor\`eme des fonctions implicites. Citons les exemples fondamentaux suivants~: les espaces des phases $(T^*L,d\p\wedge d\q)$ de la m\'ecanique classique, associ\'es \`a des espaces de configurations $L$ qui sont des vari\'et\'es lisses, munis de la structure r\'eelle canonique $c_L:(\p,\q)\mapsto(-\p,\q)$~; l'espace projectif $\P^n$ muni de la forme de Fubini-Study et de la structure r\'eelle $conj$ donn\'ee par la conjugaison complexe~; les vari\'et\'es projectives lisses d\'efinies par des polyn\^omes homog\`enes \`a coefficients r\'eels, avec la structure induite par celle de $\P^n$. On a un mod\`ele local pour $(X,\omega,c_X)$ au voisinage de $L=\R X$~: un voisinage de $L$ est isomorphe \`a un voisinage de la section nulle dans $(T^*L,d\p\wedge d\q,c_L)$~\cite{Weinstein-Lag}.

Depuis Gromov~\cite{Gromov} nous savons qu'il est utile de regarder les vari\'et\'es symplectiques comme des analogues ``flexibles" des vari\'et\'es k\"ahleriennes. De fa\c{c}on plus pr\'ecise, soit $\cJ_\omega$ l'espace des structures presque-complexes $J$ (de classe $C^\ell$, $\ell\gg 1$) sur $TX$ qui sont 
\emph{$\omega$-compatibles}, i.e. telles que $\omega(\cdot,J\cdot)$ est une m\'etrique riemannienne. Les \'el\'ements de $\cJ_\omega$ sont les sections $C^\ell$ d'un fibr\'e \`a fibres contractiles, isomorphes \`a $\mathrm{Sp}(2n)/\mathrm{U}(n)$, de sorte que $\cJ_\omega$ est une vari\'et\'e de Banach s\'eparable non-vide et contractile. En pr\'esence d'une structure r\'eelle $c_X$ on obtient une involution 
$
\overline{c_X}^*:\cJ_\omega\to\cJ_\omega,\quad J\mapsto -dc_X\circ J\circ dc_X 
$
dont le lieu des points fixes $\R \cJ_\omega$ est l'espace des structures presque-complexes $\omega$-compatibles qui rendent $c_X$ anti-holomorphe. Welschinger a montr\'e que l'espace $\R \cJ_\omega$ est une vari\'et\'e de Banach s\'eparable et contractile~\cite[\S1.1]{We-Inv}.

Un choix de $J\in\cJ_\omega$ permet de consid\'erer l'espace des \emph{courbes (rationnelles) $J$-holomorphes} $u:\P^1\to (X,J)$, solutions de l'\'equation 
$du+J\circ du\circ j=0$ o\`u $j$ est la structure complexe de $\P^1$. C'est une \'equation de type Cauchy-Riemann, elliptique d'indice $2c_1(X)d+2n$, o\`u $d=u_*[\P^1]\in H_2(X;\Z)$ et $c_1(X)$ est la premi\`ere classe de Chern du fibr\'e complexe $(TX,J)$. Lorsque $J\in\R\cJ_\omega$, il existe une involution naturelle $u\mapsto c_X\circ u\circ conj$ sur l'espace des courbes $J$-holomorphes, et ses points fixes s'appellent \emph{courbes $J$-holomorphes r\'eelles}. 
On consid\`ere par la suite les espaces de solutions modulo reparam\'etrisation conforme \`a la source, et on parle alors d'\emph{espaces de modules}. 

La vari\'et\'e symplectique $(X,\omega)$ est dite \emph{semi-positive} (resp. \emph{fortement semi-positive}) si, pour toute classe sph\'erique $d\in H_2(X;\Z)$ telle que $[\omega]d>0$, on a l'implication $c_1(X)d\ge 3-n \ \Rightarrow \ c_1(X)d\ge 0$ (resp. $c_1(X)d\ge 2-n \ \Rightarrow \ c_1(X)d\ge 0$). Les invariants de Gromov-Witten (en genre $0$) d'une vari\'et\'e semi-positive sont d\'efinis de la fa\c{c}on suivante~\cite{MS04}~: on se  restreint aux courbes dites \emph{simples} qui ne factorisent pas \`a travers un rev\^etement ramifi\'e non-trivial de $\P ^1$, on enrichit la source $\P^1$ de points marqu\'es mobiles, on impose des conditions d'incidence aux points marqu\'es de fa\c{c}on \`a ramener la dimension des espaces de modules \`a z\'ero, et finalement on compte les solutions. Notons les deux sp\'ecificit\'es suivantes~: (i) le r\'esultat peut \^etre interpr\'et\'e de fa\c{c}on duale comme le calcul d'une int\'egrale sur l'espace de modules de courbes avec points marqu\'es. Ce dernier porte une \emph{classe fondamentale} puisqu'il poss\`ede une \emph{compactification par des strates de codimension $\ge 2$}. Par ailleurs, les conditions d'incidence ne doivent pas n\'ecessairement \^etre ponctuelles, ce qui a des cons\'equences profondes comme par exemple l'existence du produit quantique sur $H^*(X;\Z)$~\cite{KM94,RT95,MS04}~; (ii) lorsque les conditions d'incidence sont repr\'esent\'ees par des sous-vari\'et\'es $J$-complexes, les courbes sont \emph{compt\'ees avec le m\^eme signe}. Ceci est une manifestation de la positivit\'e des intersections des objets holomorphes. 

Le cas particulier des conditions d'incidence ponctuelles est fondamental pour la suite. Soit $d\in H_2(X;\Z)$ et $\cM^d_k(X,J)^*$ l'espace des modules de courbes $J$-holomorphes simples avec $k$ points marqu\'es, qui repr\'esentent la classe $d$. 

\begin{center}
\emph{On suppose que $(n-1)$ divise $(c_1(X)d-2)$ et on pose $k=k_d=\frac 1{n-1}(c_1(X)d-2)+1$.} 
\end{center}

Soit $\ux\in X^k$ un $k$-uplet de points deux-\`a-deux distincts. Pour un choix g\'en\'erique de $J\in\cJ_\omega$ l'espace $\cM^d_k(X,J)^*$ est une vari\'et\'e de dimension $2c_1(X)d+2n-6+2k$ et 
$\ux$ est une valeur r\'eguli\`ere de l'application d'\'evaluation 
$$
\ev_J:\cM^d_k(X,J)^*\to X^k, \qquad (u,z_1,\dots,z_k)\mapsto (u(z_1),\dots,u(z_k)).
$$
La valeur de $k$ a \'et\'e choisie telle que la source et le but de $\ev_J$ soient de m\^eme dimension. Notons que $\ev_J$ admet une extension naturelle $\oev_J:\overline\cM^d_k(X,J)^*\to X^k$ \`a la compactification de Gromov-Kontsevich par des courbes stables~\cite{KM94,MS04}. On note $\cM^d(\ux,J)=\ev_J^{-1}(\ux)$ et $\overline\cM^d(\ux,J)=\oev_J^{-1}(\ux)$. Pour $J$ g\'en\'erique la fibre $\cM^d(\ux,J)$ est compacte, form\'ee d'un nombre fini de points qui sont des courbes immerg\'ees. Voici une \'ebauche d'argument pour montrer que ce nombre $N_d(\underline x, J)$ ne d\'epend pas du choix de $\underline x$ ou du choix g\'en\'erique de $J$. On fixe $\underline x$ et on d\'emontre l'ind\'ependance par rapport \`a $J$ en utilisant une m\'ethode  de continuit\'e~: soit $\cJ_\omega^0(\ux)\subset\cJ_\omega$ l'ensemble des $J$ pour lesquels $\cM^d(\ux,J)$ contient des courbes cuspidales, ou pour lesquels $\overline\cM^d(\ux,J)$ contient des courbes r\'eductibles. Le point cl\'e est que $\cJ_\omega^0(\ux)$ est une union au plus d\'enombrable de sous-vari\'et\'es de codimension $\ge 2$ de $\cJ_\omega$. Soient $J_0,J_1\in\cJ_\omega$ tels que $N_d(\ux, J_0)$ et $N_d(\ux, J_1)$ soient d\'efinis. L'espace $\cJ_\omega$ \'etant contractile, il existe un chemin $J_t$, $t\in[0,1]$ reliant $J_0$ \`a $J_1$. Pour un choix g\'en\'erique du chemin l'espace de modules \`a param\`etre 
$\cM=\bigcup_t\, \{t\}\times \cM^d(\ux,J_t)$ est une vari\'et\'e de dimension $1$, et les points critiques de la projection naturelle $\cM\to[0,1]$ sont les courbes cuspidales. Un chemin g\'en\'erique \'evite $\cJ_\omega^0(\ux)$, de sorte que $\cM$ est un rev\^etement propre de $[0,1]$. On obtient que la fonction $N_d(\ux, J_t)$ est localement constante, donc constante sur $[0,1]$. Par ailleurs $N_d(\ux,J)$ est localement constante en $\underline x$ \`a $J$ fix\'e. Puisque $X^k\setminus \mathrm{Diag}$ est connexe, on obtient que $N_d=N_d(\underline x,J)$ ne d\'epend pas du choix de $\underline x$ et $J$. 

Supposons maintenant donn\'ee une structure r\'eelle $c_X$ et une classe d'homologie $d\in H_2(X;\Z)$ telle que $(c_X)_*d=-d$ et $(n-1)$ divise $c_1(X)d-2$, de sorte que $k=k_d\in\N$. On suppose par la suite $k\ge 1$. On consid\`ere $J\in\R\cJ_\omega$ et une collection r\'eelle de points $\ux\in X^k\setminus\mathrm{Diag}$, i.e. une collection compos\'ee de $r$ points dans $\R X$ et de $r_X$ paires de points dans $X\setminus \R X$ conjugu\'es par $c_X$, avec $r+2r_X=k$. Notons $R^d(\ux,J)$ le nombre de courbes rationnelles $J$-holomorphes r\'eelles passant par $\ux$. 

Le nombre $R^d(\ux,J)$ d\'epend \`a $r$ fix\'e du choix de $\ux$ ou encore, de fa\c{c}on \'equivalente, du choix de $J$\footnote{Voir l'exemple des cubiques rationnelles dans $\P^2$ \`a la fin de cette section.}. Il s'ensuit que l'argument qui montrait l'ind\'ependance de $N_d(\ux,J)$ par rapport aux choix doit n\'ecessairement tomber en d\'efaut. En effet, ce dernier \'etait bas\'e sur le fait que l'espace $\cJ^0_\omega(\ux)$ des ``accidents" vivait en codimension $\ge 2$. Par contraste, la partie r\'eelle $\R\cJ^0_\omega\subset \R\cJ_\omega$ vit en codimension $\ge 1$ et ne pourra plus \^etre \'evit\'ee par un chemin $(J_t)\subset \R\cJ_\omega$ g\'en\'erique. 

Ceci \'etait \`a peu de choses pr\`es la situation avant les travaux 
que nous exposons dans cet article. Il y avait des murs -- dans $\R\cJ_\omega$ -- que l'on ne savait pas comment franchir. Welschinger a imagin\'e le ph\'enom\`ene suivant~: 

\smallskip 

\begin{center}
\emph{Il est possible d'attribuer des signes aux courbes r\'eelles soumises \`a des conditions d'incidence ponctuelles de mani\`ere \`a ce que leur comptage alg\'ebrique soit un invariant.} 
\end{center}

\smallskip

%

La d\'emarche de Welschinger est la suivante~: (i)~il fixe une composante connexe $L$ de $\R X$ et un entier $0\le r\le k_d$, avec $r\ge 1$ si $n\ge 3$~; (ii) il fixe une collection r\'eelle de points $\ux\in\R(X^k\setminus\mathrm{Diag})$ ayant $r$ points r\'eels appartenant tous \`a $L$, et il choisit $J\in\R\cJ_\omega$ assez g\'en\'erique pour que $\ux$ soit une valeur r\'eguli\`ere de l'application d'\'evaluation 
$$
\R\ev:\R\cM^d_k(X,J)^*\to\R(X^k)
$$ 
et pour que la pr\'eimage de $\ux$ par $\R\oev:\R\overline\cM^d_k(X,J)^*\to\R(X^k)$ ne contienne pas de courbes r\'eductibles. Ceci assure en particulier la compacit\'e de $\R\ev^{-1}(\ux)$~; (iii) \emph{il attribue des signes} $\eps^\pp(C)\in\{\pm1\}$ \emph{aux \'el\'ements} $C\in\R\ev^{-1}(\ux)$ et d\'efinit 
$$
\chi^{d,\pp}_r(\ux,J):=\sum_{C\in\R\ev^{-1}(\ux)} \eps^\pp(C).
$$
On appellera les signes $\eps^\pp(C)$ \emph{signes de Welschinger\footnote{Lorsque $n\ge 3$, Welschinger appelle ce signe \emph{\'etat spinoriel de $C$} et le note $sp(C)$~\cite{We-Duke,We-semipositive}.}}. La recette d'attribution des signes est valable pour des vari\'et\'es symplectiques de dimension arbitraire. En dimension $n\ge 3$, on requiert qu'une  certaine hypoth\`ese de nature topologique soit v\'erifi\'ee, hypoth\`ese qui assure l'existence d'une structure $\mathrm{Pin}_n^-$, ou $\mathrm{Pin}_n^+$, ou $\mathrm{Spin}$ sur un fibr\'e vectoriel appropri\'e au-dessus de $\R X$. La d\'efinition des signes en dimension $n\ge 3$ utilise le choix d'une telle structure $\pp$, cf.~\S\ref{sec:dim3}. Finalement, (iv)~Welschinger d\'emontre 

\begin{theo} [\cite{We-Inv,We-Duke,We-semipositive}] \label{thm:main} 
Soit $L$ une composante connexe de $\R X$ et $r\in\{0,\dots,k_d\}$ avec $r\ge 1$ si $n\ge 3$. On consid\`ere des collections $\ux\in\R(X^k\setminus\mathrm{Diag})$ ayant $r$ points r\'eels, situ\'es tous sur $L$. En dimension $2$ (resp. $3$), le nombre $\chi^d_r(L)=\chi^d_r(\ux,J)$ (resp. $\chi^{d,\pp}_r(L)=\chi^{d,\pp}_r(\ux,J)$ avec $\pp$ une structure $\mathrm{Pin}_3^-$) ne d\'epend ni du choix de $\ux$, ni du choix g\'en\'erique de $J$.
\end{theo}

Puisque le lieu r\'eel de $\P^1$ est connexe, il ne peut y avoir de courbe $J$-holomorphe r\'eelle passant par $\ux$ que si tous les points r\'eels de $\ux$ appartiennent \`a la m\^eme composante de $\R X$. Ceci justifie le choix d'une composante $L$ dans la d\'efinition de l'invariant $\chi^d_r(L)$, qui d\'epend par ailleurs de ce choix. Lorsque $\R X$ est connexe on note simplement $\chi^d_r(\R X)=\chi^d_r$. L'estim\'ee fondamentale suivante d\'ecoule directement de la d\'efinition~:

\begin{coro} \label{cor:borneinf} Le nombre $R^d(\ux,J)$ de courbes $J$-holomorphes r\'eelles passant par $\ux$ v\'erifie 
$$
|\chi^d_r(L)|\le R^d(\underline x,J)\le N_d.
$$
\end{coro}

\smallskip 

\smallskip

\noindent \emph{Exemple. --- Cubiques rationnelles dans $\P^2$}~\cite[Prop.~4.7.3]{DK}, cf. aussi \cite[Prop.~3.6]{Sottile}.
Soit $X=\P^2$ et $d=3[\P^1]$, de sorte que $k_d=8$, $r\in\{0,2,4,6,8\}$ et $r_X\in\{4,3,2,1,0\}$.

(i) Le nombre de cubiques rationnelles \emph{complexes} passant par $8$ points g\'en\'eriques de $\P^2$ vaut $12$. C'est un calcul classique qui remonte \`a Schubert. On consid\`ere un pinceau de cubiques dont le lieu de base contient ces points et on note $Z$ l'\'eclatement de $\P^2$ aux $9$ points du lieu de base, de sorte que $\chi(Z)=12$. Pour un choix g\'en\'erique des points les fibres singuli\`eres de la fibration $Z\to\P^1$ sont nodales ($\chi=1$), alors que les fibres g\'en\'eriques sont des courbes elliptiques ($\chi=0$). Ce ne sont donc que les fibres singuli\`eres qui contribuent \`a la caract\'eristique d'Euler de $Z$ et il doit y en avoir exactement $12$.  

(ii) Le nombre $R^3(\ux)$ de cubiques rationnelles \emph{r\'eelles} passant par une collection $\ux$ g\'en\'erique compos\'ee de $r$ points r\'eels et $r_X$ paires de points complexes conjugu\'es est minor\'e par $r$. Pour le voir, consid\'erons un pinceau de cubiques r\'eelles dont le lieu de base est constitu\'e de ces points et d'un point additionnel qui est r\'eel. La partie r\'eelle $\R Z$ de l'\'eclatement fibre au-dessus de $\R\P^1$, et on a $\chi(\R Z)=\chi(\R\P^2)-(r+1)=-r$.
Les fibres lisses sont des parties r\'eelles de cubiques lisses, hom\'eomorphes \`a une ou deux copies du cercle ($\chi=0$). Ce ne sont donc que les fibres singuli\`eres qui contribuent \`a la caract\'eristique d'Euler. Celles-ci sont des cubiques r\'eelles nodales, ayant un unique n\oe ud r\'eel, et leur caract\'eristique d'Euler vaut $\pm 1$ selon que ce n\oe ud r\'eel est solitaire (intersection de deux branches complexes conjugu\'ees) ou non-solitaire (intersection de deux branches r\'eelles). Dans le premier cas la cubique est l'union d'un point et d'un cercle, dans le deuxi\`eme cas la cubique est hom\'eomorphe \`a un bouquet de deux cercles. Notons $c_-$ le nombre de fibres singuli\`eres avec un n\oe ud solitaire et $c_+$ le nombre de fibres singuli\`eres avec un n\oe ud non-solitaire, de sorte que 
$$
c_- - c_+=-r.
$$ 
Alors $R^3(\ux)=c_- + c_+$ v\'erifie l'estim\'ee 
$$
r\le R^3(\ux)\le 12. 
$$
Pour $r=8$ chacune des valeurs possibles $8,10,12$ est atteinte pour un choix appropri\'e de $\ux$~\cite{DK,Sottile}. 
Les signes de Welschinger pour les cubiques rationnelles r\'eelles sont 
$$
\eps(C)=\left\{\begin{array}{ll} 
-1, & \mbox{si } C \mbox{ a un n\oe ud r\'eel solitaire,} \\
+1, & \mbox{si } C  \mbox{ a un n\oe ud r\'eel non-solitaire.}
\end{array}\right.
$$
La valeur de l'invariant de Welschinger est donc
$$
\chi^3_r = c_+ - c_- =  r = -\chi(\R Z).
$$

\smallskip 

\noindent {\it Remarque.}{\ }--- Cette identification avec la caract\'eristique d'Euler du lieu r\'eel d'une vari\'et\'e, bien qu'amusante puisque les notations co\"{\i}ncident, semble fortuite. La notation $\chi^d_r$ utilis\'ee par Welschinger est plut\^ot motiv\'ee par l'espoir que cet invariant puisse \^etre ``cat\'egorifi\'e" en l'interpr\'etant comme caract\'eristique d'Euler d'un groupe d'homologie. La topologie des petites dimensions regorge de tels exemples, dont le plus fameux est celui de l'homologie de Khovanov qui cat\'egorifie le polyn\^ome de Jones. 

\smallskip

Notre article est structur\'e de la mani\`ere suivante. La section \S\ref{sec:chird} contient la d\'efinition des signes de Welschinger en dimension arbitraire et la preuve du th\'eor\`eme~\ref{thm:main}. Dans \S\ref{sec:optimalite} on explore une deuxi\`eme direction de recherche initi\'ee par Welschinger qui est celle de l'utilisation des techniques de la th\'eorie symplectique des champs~\cite{BEHWZ,EGH}  en g\'eom\'etrie alg\'ebrique r\'eelle\footnote{Nous nous devons de mentionner un th\'eor\`eme pr\'ecurseur d\^u \`a Viterbo~\cite{V} et Eliashberg~\cite{EGH}, expliqu\'e par Kharlamov dans un expos\'e au s\'eminaire Bourbaki~\cite{slava-claude}.}. Nous nous concentrons sur un r\'esultat d'optimalit\'e pour les invariants $\chi^d_r$ en dimension quatre et nous donnons quelques ouvertures vers les invariants relatifs.  
Finalement, dans \S\ref{sec:developpements} nous donnons un aper\c{c}u de r\'esultats connexes aux invariants de Welschinger, obtenus notamment en utilisant des techniques de g\'eom\'etrie tropicale ou des id\'ees inspir\'ees par la conjecture de sym\'etrie miroir. 

Le lecteur est chaleureusement invit\'e \`a consulter l'article de survol~\cite{We10} \'ecrit par Welschinger \`a l'occasion de l'ICM 2010.

\section{Invariants pour des conditions d'incidence ponctuelles} \label{sec:chird}

\subsection{L'invariant de Welschinger en dimension $2$} On consid\`ere une vari\'et\'e symplectique r\'eelle de dimension quatre $(X,\omega,c_X)$ et une classe d'homologie $d\in H_2(X;\Z)$ telle que $(c_X)_*d=-d$. (Les vari\'et\'es de dimension quatre sont automatiquement fortement semi-positives.) On note $k=k_d=c_1(X)d-1$ et l'on suppose $k_d\in\N^*$.  On fixe une collection r\'eelle $\ux\in X^k\setminus\mathrm{Diag}$, compos\'ee de $r$ points r\'eels situ\'es sur la m\^eme composante de $\R X$ et de $r_X$ paires de points dans $X\setminus \R X$ conjugu\'es par $c_X$, avec $r+2r_X=k$. 

Soit $J\in\R\cJ_\omega$. On note $\cR^d(\ux,J)$ l'ensemble des courbes $J$-holomorphes simples homologues \`a $d$ passant par $\ux$. Pour $J$ g\'en\'erique, cet ensemble est fini et consiste en des courbes irr\'eductibles, immerg\'ees, ayant des points doubles transverses (n\oe uds). (Par la formule d'adjonction, il y en a $\frac 1 2 (d^2-c_1(X)d+2)$.)
Un n\oe ud r\'eel peut \^etre solitaire (intersection de deux branches complexes conjugu\'ees), ou pas (intersection de deux branches r\'eelles). Welschinger d\'efinit la \emph{masse $m(C)$ d'une courbe $C\in\cR^d(\ux,J)$} comme le nombre de n\oe uds r\'eels solitaires. Le \emph{signe de Welschinger de} $C$ est d\'efini par
$$
\eps(C)=(-1)^{m(C)}.
$$
Le lecteur est invit\'e \`a comparer cette d\'efinition avec l'exemple de la section pr\'ec\'edente. On pose 
$$
\chi^d_r(\ux,J)=\sum_{C\in\cR^d(\ux,J)} \eps(C).
$$

\noindent {\it D\'emonstration du th\'eor\`eme~\ref{thm:main} dans le cas $n=2$~\cite{We03,We-Inv}.} --- L'ind\'ependance de $\chi^d _r(\underline x,J)$ par rapport au choix de $\underline x$ est une cons\'equence de l'ind\'ependance par rapport au choix de $J$ \`a $\underline x$ fix\'e. Ceci d\'ecoule de l'observation suivante~: deux collections $\underline x$, $\underline x'$ dont les points r\'eels sont en nombre \'egal et sont situ\'es sur la m\^eme composante connexe de $\R X$ sont reli\'ees par une isotopie symplectique r\'eelle, qui induit en particulier une isotopie de structures presque complexes r\'eelles. 

Fixons maintenant une collection r\'eelle $\underline x$ et deux structures presque complexes r\'eelles g\'en\'eriques $J_0$ et $J_1$. L'espace des structures r\'eelles compatibles avec $\omega$ \'etant contractile, il existe un chemin lisse $J:[0,1]\to\cJ^\ell_\omega$ reliant $J_0$ \`a $J_1$. Pour un choix g\'en\'erique du chemin $(J_t)$ l'espace de modules \emph{\`a param\`etre} 
$$
\R \cM:=\bigsqcup_{t\in[0,1]} \{t\}\times \cR^d(\underline x,J_t)
$$ 
est une vari\'et\'e lisse de dimension $1$. La premi\`ere projection $\pi:\R\cM\to[0,1]$ est une application de Fredholm d'indice $0$. 
Welschinger d\'emontre dans \cite{We-Inv} que, lorsque le chemin $(J_t)$ est choisi de fa\c{c}on g\'en\'erique, il existe un ensemble fini $0<t_1<\dots<t_N<1$ qui v\'erifie les conditions suivantes~: 
\begin{enumerate}
\item le nombre $\chi^d_r(\underline x,J_t)$ est d\'efini et est localement constant pour $t\neq t_i$~; 
\item pour $t=t_i$ l'une des situations suivantes se produit~: 
 \begin{enumerate}
\item $\cR^d(\underline x,J_{t_i})$ contient une courbe dont toutes les singularit\'es sont des points doubles ordinaires \`a l'exception de l'une d'entre elles qui est un point de tangence de deux branches, ou un point triple r\'eel ordinaire, ou un point de rebroussement r\'eel de premi\`ere esp\`ece~; 
\item il existe une suite $t_\nu\to t_i$ et des \'el\'ements $C_\nu\in\cR^d(\underline x,J_{t_\nu})$ qui convergent au sens de Gromov vers une courbe r\'eductible compos\'ee de deux branches irr\'eductibles r\'eelles dont les singularit\'es sont des points doubles ordinaires~;
 \end{enumerate}
\end{enumerate} 
(Intuitivement, la convergence de Gromov g\'en\'eralise au cadre $J$-holomorphe la d\'eg\'en\'erescence d'une conique de $\P^2$ vers une union de deux droites.) L'apparition du cas~(b) est \'equivalente \`a la non-compacit\'e de $\R\cM$, alors que les courbes cuspidales du~(a) correspondent aux points critiques de la projection $\pi$. Ceux-ci peuvent \^etre suppos\'es non-d\'eg\'en\'er\'es pour un choix g\'en\'erique du chemin $(J_t)$~; ce seront en particulier des maxima ou des minima locaux pour $\pi$. 
 
Il s'agit de montrer que l'invariant $\chi^d_r(\underline x,J_t)$ ne change pas lorsque l'on traverse une valeur $t_i$. La preuve est contenue dans la Figure~1 lorsqu'il s'agit d'un point de tangence de deux branches (la masse de la courbe en question change de $2$, $-2$, ou $0$), ou bien lorsqu'il s'agit d'un point triple r\'eel ordinaire (la masse de la courbe reste constante). 

\begin{figure}[h]
\label{fig:easy}
\centering
\includegraphics[scale=0.62]{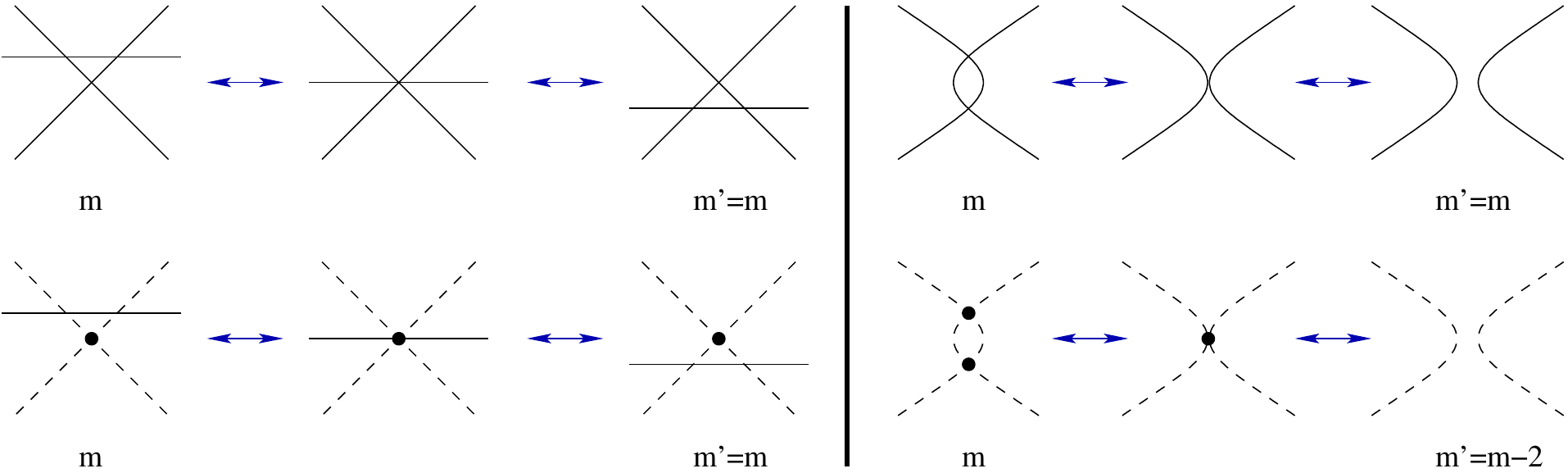}
\caption{Le mur des points triples et celui des branches tangentes.}
\end{figure}
 
Les deux autres situations sont tr\`es d\'elicates. Welschinger s'appuie dans son analyse de fa\c{c}on essentielle sur les travaux de Ivashkovich et Shevchishin~\cite{IS,Sh} pour d\'ecrire le voisinage d'une courbe r\'eductible, respectivement celui d'une courbe singuli\`ere, cf.~\S\ref{sec:courbes}. 

\noindent ($\alpha$) \emph{Le cas d'une courbe r\'eductible.} Soit $C$ une courbe $J_{t_i}$-holomorphe r\'eelle r\'eductible qui est l'union de deux composantes r\'eelles irr\'eductibles $C_1$, $C_2$ ayant $p$ points d'intersection transverse r\'eels. Welschinger d\'emontre~\cite[Proposition~2.14]{We-Inv} qu'il existe un voisinage $\cW$ de $C$ dans la compactification de Gromov $\overline \cR^d(\underline x)=\bigsqcup_J \{J\}\times \overline \cR^d(\underline x,J)$ tel que, pour $t$ proche de $t_i$, l'intersection $\cR^d(\underline x,J_t)\cap \cW$ consiste en exactement $p$ courbes $J_t$-holomorphes r\'eelles irr\'eductibles, chacune d'entre elles \'etant obtenue topologiquement en lissant un des points r\'eels d'intersection de $C_1$ et $C_2$. En fixant $t^-<t_i$ et $t^+>t_i$ proches de $t_i$ on obtient de cette mani\`ere deux collections compos\'ees de $p$ courbes chacune, qui sont naturellement mises en bijection. Les courbes qui se correspondent ont le m\^eme nombre de points doubles r\'eels solitaires, de sorte que leurs masses sont \'egales. Par cons\'equent le nombre $\chi^d_r(\underline x,J_t)$ ne change  pas lorsque l'on traverse la valeur $t_i$.

\noindent ($\beta$) \emph{Le cas d'une courbe cuspidale.} Soit $C\in \cR^d(\underline x, J_{t_i})$ une courbe r\'eelle ayant un unique point de rebroussement r\'eel, correspondant \`a un maximum local (resp. minimum local) non-d\'eg\'en\'er\'e de $\pi$. Welschinger d\'emontre~\cite[Proposition~2.16]{We-Inv} qu'il existe un voisinage $\cW$ de $C$ dans $\overline \cR^d(\underline x)$ tel que, pour $t<t_i$ (resp. $t>t_i$) proche de $t_i$, l'intersection $\cR^d(\underline x,J_t)\cap \cW$ consiste en exactement deux courbes $C^+$, $C^-$ telles que 
$m(C^+)=m(C^-)+1$ et, pour $t>t_i$ (resp. $t<t_i$) proche de $t_i$, l'intersection $\cR^d(\underline x,J_t)\cap \cW$ est vide, cf. Figure~2 ci-dessous. (Dans cet \'enonc\'e, la partie difficile est de montrer la relation entre les masses des courbes $C_\pm$, qui est expliqu\'ee par la remarque ci-dessous.) On en d\'eduit que le nombre $\chi^d_r(\underline x,J_t)$ ne change  pas lorsque l'on traverse la valeur $t_i$. Le fait d'avoir d\'efini $\chi^d_r$ par une somme altern\'ee s'av\`ere crucial pour cette \'etape de la preuve.
\hfill{$\square$}

\noindent {\it Remarque.} --- L'exemple suivant est fondamental pour comprendre le cas ($\beta$) ci-dessus. Consid\'erons la famille $u_\lambda:\C\to\C^2$, $z\mapsto (z^2,z^3+\lambda z)$ de courbes r\'eelles affines, index\'ee par $\lambda$ r\'eel proche de z\'ero. La courbe $u_0$ a un point de rebroussement en $z=0$, alors que $u_\lambda$ a un n\oe ud r\'eel non-solitaire (resp. solitaire) pour $\lambda<0$ (resp. $\lambda>0$) correspondant aux param\`etres $z=\pm\sqrt{-\lambda}$ (Figure~2). Lorsque l'on regarde ces courbes dans $\P^2$ leurs masses v\'erifient donc la relation $m(u_\lambda)=m(u_{-\lambda})+1$ pour $\lambda>0$. Le cas ($\beta$) d\'ecoule du fait que cette famille de courbes constitue un mod\`ele local pour $\R\cM$ au voisinage de $C\in \cR^d(\underline x, J_{t_i})$ (\cite[Lemme~2.6]{We-Inv}, \cite[Corollaire~1.4.3]{IS}).

\begin{figure}[h]
\label{fig:cusp}
\centering
\begin{tabular}{lcr}
\includegraphics[scale=0.7]{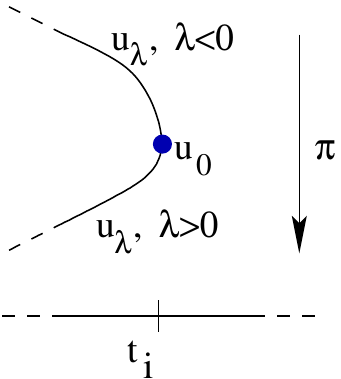}
& \qquad {\ } & 
\includegraphics[scale=0.7]{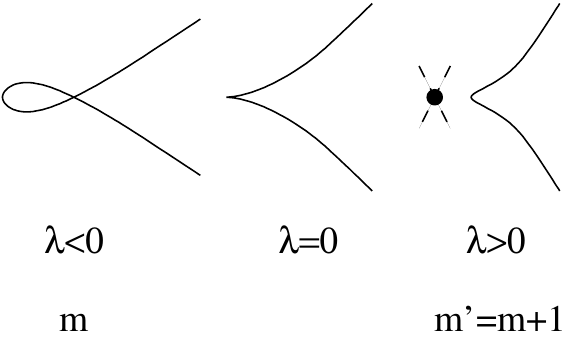}
\end{tabular}
\caption{Mod\`ele local pour le mur des courbes cuspidales.}
\end{figure}

\subsection{Espaces de modules de courbes $J$-holomorphes} \label{sec:courbes}

Nous venons de voir que la d\'efinition des signes de Welschinger en dimension $2$ est topologique. Par contraste, leur d\'efinition en dimension $n\ge 3$ m\'elange analyse et g\'eom\'etrie~: on utilise l'op\'erateur de Fredholm obtenu en lin\'earisant l'\'equation des courbes $J$-holomorphes, mais aussi une structure suppl\'ementaire de nature spinorielle sur $\R X$. Dans cette section nous introduisons les notions analytiques fondamentales requises pour d\'efinir les signes de Welschinger en suivant les r\'ef\'erences~\cite[\S3 et App.~C]{MS04}, \cite[\S1]{IS}, \cite[\S1]{Barraud} et~\cite{We-Duke}. On choisit $d\in H_2(X;\Z)$ comme avant et on suppose $k=k_d\in\N^*$. 

Fixons $p>2$. Soit $\cB=\cB^d$ la vari\'et\'e de Banach dont les points sont les applications $u:\P^1\to X$ de classe de Sobolev $W^{1,p}$ homologues \`a $d$. Les \'el\'ements de $\cB$ sont en particulier des fonctions continues. Un choix de $J\in\cJ_\omega$ d\'etermine un fibr\'e de Banach $\cE\to \cB$ de fibre $\cE_u=L^p(\Lambda^{0,1}\P^1\otimes u^*TX)$, l'espace des formes $(j,J)$-antilin\'eaires de classe $L^p$ \`a valeurs dans $u^*TX$. Le lieu des z\'eros de la section
$$
\dbar_J=\frac  1 2 (du+J(u)\circ du\circ j):\cB\to \cE,
$$
not\'e $\widehat \cM^d(X,J)$, est l'espace des applications $J$-holomorphes lisses homologues \`a $d$, la lissit\'e \'etant une cons\'equence de la r\'egularit\'e elliptique. On note $\widehat \cM^d(X,J)^*\subset \widehat \cM^d(X,J)$ le sous-espace des courbes $J$-holomorphes simples. Pour un choix g\'en\'erique de $J\in\cJ_\omega$ la section $\dbar_J$ est transverse \`a la section nulle le long de $\widehat \cM^d(X,J)^*$~; un tel $J$ est dit \emph{r\'egulier}. La r\'egularit\'e de $J$  \'equivaut \`a la surjectivit\'e en tout point $u\in\widehat \cM^d(X,J)^*$ de la composante verticale de la diff\'erentielle $d\dbar_J(u)^{\mathrm{vert}}:T_u\cB\to \cE_u$, not\'ee
$$
D_u=D_{u,J}:W^{1,p}(u^*TX)\to L^p(\Lambda^{0,1}\P^1\otimes u^*TX).
$$
Choisissons des coordonn\'ees locales sur $X$ et $\P^1$, conformes sur ce dernier. La diff\'erentielle verticale s'\'ecrit alors sous la forme $D_u\xi=\dbar_J\xi - \frac 1 2 (J \partial_\xi J)(u)\partial _J(u)$. On voit en particulier que $D_u$ est un op\'erateur de Cauchy-Riemann \emph{g\'en\'eralis\'e}, au sens o\`u il v\'erifie la relation $D_u(f\xi)=fD_u\xi+\dbar f\otimes\xi$ pour toute fonction $f$ sur $\P^1$ \`a valeurs \emph{r\'eelles} \footnote{Les op\'erateurs de Cauchy-Riemann sont les op\'erateurs $W^{1,p}(u^*TX)\to L^p(\Lambda^{0,1}\P^1\otimes u^*TX)$ qui v\'erifient cette identit\'e pour toute fonction $f$ \`a valeurs complexes.}. L'op\'erateur $D_u$ est en particulier elliptique et son indice (r\'eel) est donn\'e par la formule de Riemann-Roch
$$
\mathrm{ind} \, D_u = 2n+2c_1(X)d.
$$
La dimension de $\widehat \cM^d(X,J)^*$ est \'egale \`a $\mathrm{ind}\, D_u$ lorsque $D_u$ est surjectif. Le groupe $\mathrm{PGL}(2,\C)$ agit librement sur $\widehat \cM^d(X,J)$ par reparam\'etrisation \`a la source, et l'espace de modules de courbes $J$-holomorphes  
$
\cM^d(X,J)^*:=\widehat \cM^d(X,J)^*/\mathrm{PGL}(2,\C)
$ 
est de dimension $\mathrm{ind}\, D_u - 6= 2c_1(X)d +2n-6$. 

\smallskip 

\noindent {\it Remarque (courbes simples).} --- Nous avons restreint notre attention aux courbes simples, pour lesquelles les $J$ r\'eguliers sont g\'en\'eriques. Ceci ne r\'esulte pas en une perte de g\'en\'eralit\'e tant que l'objet d'int\'er\^et est l'image des courbes dans $X$. En effet, toute courbe $J$-holomorphe est rev\^etement ramifi\'e d'une courbe simple~\cite[\S2.5]{MS04}. 

\smallskip 

\noindent {\it Remarque (g\'en\'ericit\'e).} --- La g\'en\'ericit\'e des $J$ r\'eguliers pour les courbes simples d\'ecoule de l'argument suivant. Soit $\cB^*\subset \cB$ le sous-espace des applications dites \emph{quelque part injectives}, pour lesquelles il existe $z\in\P^1$ tel que $u^{-1}(u(z))=\{z\}$. Consid\'erons le fibr\'e de Banach $\cE\to \cB^*\times \cJ_\omega$ de fibre $\cE_{(u,J)}=L^p(\Lambda^{0,1}\P^1\otimes u^*TX)$ muni de la section $\dbar(u,J)=\dbar_J(u)$. Le point cl\'e est que la section $\dbar$ est transverse \`a la section nulle, de sorte que l'espace de modules universel $\cP^*=\dbar^{-1}(0)\subset \cB^*\times \cJ_\omega$ est une sous-vari\'et\'e de Banach. La projection $\pi:\cP^*\to\cJ_\omega$ sur le deux\`eme facteur v\'erifie $\ker d\pi_{(u,J)}\simeq \ker D_{u,J}$, $\coker d\pi_{(u,J)}\simeq \coker D_{u,J}$, de sorte que $\pi$ est une application de Fredholm de m\^eme indice que $D_{u,J}$. Le th\'eor\`eme de Sard-Smale assure que les valeurs r\'eguli\`eres $J$ de $\pi$ forment un ensemble dense, et ces $J$ sont en particulier r\'eguliers au sens pr\'ec\'edent. \emph{Tous les \'enonc\'es de g\'en\'ericit\'e de cet article se d\'emontrent selon un sch\'ema similaire, en construisant un espace de modules universel appropri\'e et en calculant l'indice de Fredholm de la projection $\pi$.} \`A titre d'exemple, mentionnons le fait que les \'el\'ements de $\cR^d(\ux,J)$ sont immerg\'es lorsque $J$ est g\'en\'erique, ou bien le fait qu'un chemin $(J_t)$ g\'en\'erique rencontre des murs d'un certain type. 

\smallskip 

On a l'identit\'e remarquable $D_u\circ du =du\circ \dbar$~\cite[Lemme~1.3.1]{IS}. Celle-ci d\'etermine un diagramme commutatif 
$$
\xymatrix{0\ar[r]Ê& W^{1,p}(T\P^1)\ar[r]^-{du} \ar[d]^{\dbar} & W^{1,p}(E) \ar[r] \ar[d]^{D_u} & W^{1,p}(E)/\mathrm{im}\,du \ar[r]Ê\ar[d]^{\oD_u} & 0\\
0\ar[r] & L^p(\Lambda^{0,1}\P^1)\ar[r]^-{du} & L^p(\Lambda^{0,1}\P^1\otimes E) \ar[r] & L^p(\Lambda^{0,1}\P^1\otimes E)/\mathrm{im}\,du \ar[r] & 0
}
$$
L'op\'erateur induit $\oD_u$ est de Fredholm d'indice $\ind\,\oD_u=\ind\,D_u-\ind\,\dbar=\ind\,D_u-6$. Si $D_u$ est surjectif alors $\oD_u$ l'est aussi et l'espace tangent \`a $\cM^d(X,J)^*$ s'identifie naturellement \`a $\ker\,\oD_u$. 

En consid\'erant la partie $\C$-lin\'eaire de l'op\'erateur $D_u$ on aboutit \`a une description fine du noyau et du conoyau de $\oD_u$. On note $\dbar_{u,J}$ (resp. $R$) la partie $\C$-lin\'eaire (resp. $\C$-anti-lin\'eaire) de $D_u$. L'op\'erateur $R$ est d'ordre z\'ero\footnote{Plus pr\'ecis\'ement $R\xi=\frac 1 4 N_J(\xi,du(\cdot))$, avec $N_J(X,Y)=[X,Y]+J[JX,Y]+J[X,JY]-[JX,JY]$ le tenseur de Nijenhuis~\cite[Lemme~C.7.3]{MS04}.}. L'op\'erateur de Cauchy-Riemann $\dbar_{u,J}:W^{1,p}(u^*TX)\to L^p(\Lambda^{0,1}\P^1\otimes u^*TX)$ d\'efinit une structure holomorphe unique sur $u^*TX$ dont c'est l'op\'erateur $\dbar$ canonique~\cite[Lemme~1.2.3]{IS} (voir aussi~\cite[\S{I}.3]{Kobayashi}). On note $E=u^*TX$ le fibr\'e holomorphe d\'efini par $\dbar_{u,J}$. L'identit\'e $D_u\circ du =du\circ \dbar$ implique $\dbar_{u,J}\circ du=du\circ \dbar$ puisque $du$ est $\C$-lin\'eaire. Lorsque $du\not\equiv 0$ on en d\'eduit un morphisme analytique injectif $\cO(T\P^1)\stackrel{du}\longrightarrow \cO(E)$ qui s'ins\`ere dans une suite exacte courte 
$$
0\longrightarrow \cO(T\P^1)\stackrel {du}\longrightarrow \cO(E)\longrightarrow \cN_u\longrightarrow 0 .
$$
Le faisceau $\cN_u$ se d\'ecompose 
$$
\cN_u=\cO(N_u)\oplus \cN^{sing}_u
$$
avec $N_u$ un fibr\'e vectoriel holomorphe de rang $n-1$, que l'on appelle \emph{fibr\'e normal de $u$}, et $\cN_u^{sing}$ un faisceau gratte-ciel \`a support dans l'ensemble des z\'eros de $du$. La fibre de $\cN_u^{sing}$ en $p$ est $\C^{\mu_p}$, avec $\mu_p$ l'ordre d'annulation de $du$ en $p$. 

\begin{prop} \cite[Lemme~1.5.1]{IS} \label{prop:IS} On a 
$$
\ker \oD_u \simeq H^0(\cN_u) = H^0(N_u)\oplus H^0(\cN_u^{sing}), \qquad \coker \oD_u \simeq H^1(N_u).
$$ 
\end{prop}

Ce type de calcul permet de d\'ecrire la codimension g\'en\'erique de diff\'erents types d'accidents dans des familles de courbes $J$-holomorphes, cf. remarque pr\'ec\'edente. Par exemple, les courbes soumises \`a $k_d$ contraintes ponctuelles et qui ne sont pas immerg\'ees apparaissent en codimension complexe $n-1$ par rapport \`a $J$, respectivement en codimension $n-1$ r\'eelle par rapport \`a $J\in\R\cJ_\omega$~\cite[Lemme~5.1]{We-semipositive}. 
Lorsque $u$ est une immersion on a $\cN_u^{sing}=0$ et l'op\'erateur $\oD_u$ est un op\'erateur de Cauchy-Riemann g\'en\'eralis\'e sur $N_u$. Sa partie $\C$-lin\'eaire est l'op\'erateur $\dbar$ canonique agissant sur les sections de $N_u$.

Soit $\widehat \cM^d_k(X,J)^*$ l'espace des courbes $J$-holomorphes simples avec $k$ points marqu\'es, constitu\'e de paires $(u,\uz)$ avec $u\in\widehat \cM^d(X,J)^*$ et $\uz=\{z_1,\dots,z_k\}\subset \P^1$ une collection de $k$ points distincts. L'espace de modules $\cM^d_k(X,J)^*$ est le quotient par l'action diagonale de $\mathrm{PGL}(2,\C)$ via $\phi\cdot(u,z_1,\dots,z_n)=(u\circ \phi^{-1},\phi(z_1),\dots,\phi(z_n))$. On a 
$
\dim\,\cM^d_k(X,J)^*=2c_1(X)d + 2n +2k -6.
$
Le but et la source de l'application d'\'evaluation naturelle
$$
\ev_J: \cM^d_k(X,J)^*\to X^k,\qquad \ev[u,z_1,\dots,z_n]:=(u(z_1),\dots,u(z_n))
$$
sont de m\^eme dimension exactement lorsque $k=k_d$. 

On suppose d\'esormais que $u$ est une courbe immerg\'ee, de sorte que $\cN_u^{sing}=0$. Soit $W^{1,p}_{-\uz}(N_u)\subset  W^{1,p}(N_u)$ le sous-espace de codimension $2(n-1)k_d$ constitu\'e des sections qui s'annulent aux points de $\uz$. Notons $D_{\uz}$ la restriction de $\oD_u$ \`a $W^{1,p}_{-\uz}(N_u)$, un op\'erateur de Fredholm d'indice $\ind\,D_{\uz}=\ind\,\oD_u-2(n-1)k_d=0$. On note $N_{u,-\uz}=N_u\otimes\cO_{\P^1}(-\uz)$. C'est un fibr\'e dont les sections holomorphes s'identifient aux sections holomorphes de $N_u$ qui s'annulent aux points de $\uz$.

\begin{prop} \cite[Lemme~1.3]{We-Duke}, \cite{We-semipositive} Supposons $u$ immerg\'ee et $\oD_u$ surjectif. On a les identifications 
suivantes~:   
$$
\ker d\ev_J\big|_{u,\uz}\simeq H^0(N_{u,-\uz}), \qquad \coker d\ev_J\big|_{u,\uz}\simeq H^1(N_{u,-\uz}).
$$
\end{prop}

\begin{coro} \label{coro:balanced} Supposons $u$ immerg\'ee et $\oD_u$ surjectif. Les conditions suivantes sont \'equivalentes. 
\renewcommand{\theenumi}{\roman{enumi}}
\begin{enumerate}
\item $(u,\uz)\in\cM^d_{k_d}(X,J)^*$ est un point r\'egulier de $d\ev_J$~; 
\item l'op\'erateur $D_\uz$ est un isomorphisme~; 
\item le fibr\'e $N_u$ se d\'ecompose en une somme directe de fibr\'es holomorphes isomorphes
$$
N_u=\cO(k_d-1)^{\oplus n-1}.
$$
\end{enumerate}
\end{coro}

\noindent \emph{Preuve.} --- Pour montrer (i)$\Leftrightarrow$(ii) on utilise les identifications 
$$
\ker\,D_{\uz}\simeq H^0(N_{u,-\uz}),\qquad \coker \,D_{\uz}\simeq H^1(N_{u,-\uz}).
$$
Celles-ci sont implicitement contenues dans la preuve par Ivashkovich et Shevchishin de la Proposition~\ref{prop:IS}. Ainsi (i) et (ii) sont \'equivalentes \`a $H^1(N_{u,-\uz})=0$. 

Montrons (i)$\Rightarrow$(iii). Par un th\'eor\`eme de Grothendieck il existe un scindement $N_u=\cO(a_1)\oplus\cdots\oplus\cO(a_{n-1})$. La condition $H^1(N_{u,-\uz})=0$ est \'equivalente \`a $0=H^1(\cO(a_i-k_d))=H^0(\cO(k_d-a_i-2))$ pour tout $i$, ce qui force $a_i\ge k_d-1$. D'un autre c\^ot\'e, puisque $u$ est immerg\'ee le degr\'e de $N_u$ est $\sum_i a_i = c_1(X)d-2=(n-1)(k_d-1)$. On en d\'eduit que $a_1=\cdots=a_{n-1}=k_d-1$. Finalement, (iii)$\Rightarrow$(i) d\'ecoule du fait que $H^1(\cO(k_d-1)\otimes \cO_{\P^1}(-\uz))=H^1(\cO(-1))=0$. 
\hfill{$\square$}

\begin{defi} \label{def:equilibre}
Une courbe $u$ immerg\'ee est dite \emph{\'equilibr\'ee} si $N_u=\cO(k_d-1)^{\oplus n-1}$. 
\end{defi}

Les courbes \'equilibr\'ees jouent un r\^ole cl\'e dans la d\'efinition des signes de Welschinger en dimension $n\ge 3$ (cf. \S\ref{sec:dim3}).

\noindent \emph{Remarque (vari\'et\'es convexes).} --- Etant donn\'ee une vari\'et\'e symplectique (fortement semi-positive), les espaces de modules de courbes simples et les invariants de Gromov-Witten sont d\'efinis pour un choix \emph{g\'en\'erique}, typiquement non-int\'egrable, de la structure presque complexe, et ceci pour des raisons de transversalit\'e. Il existe n\'eanmoins une classe de vari\'et\'es complexes o\`u la transversalit\'e est automatique~: les vari\'et\'es convexes. 

Une vari\'et\'e complexe lisse $(X,J)$ est dite \emph{convexe} si $H^1(\P^1,u^*TX)=0$ pour tout morphisme $u:\P^1\to X$. La classe principale d'exemples est constitu\'ee des espaces homog\`enes $X=G/P$, avec $G$ un groupe de Lie et $P$ un sous-groupe parabolique~\cite[\S0.4]{FP}~: espaces projectifs, grassmanniennes, vari\'et\'es de drapeaux, quadriques lisses, produits de telles vari\'et\'es. En effet, dans cette situation $TX$ est engendr\'e par ses sections globales et il en est de m\^eme pour $u^*TX$. Les op\'erateurs $D_u$ et $\oD_u$ sont automatiquement surjectifs. Les espaces de modules $\cM^d(\ux,J)$ sont d\'efinis d\`es que $\ux$ est une valeur r\'eguli\`ere de $\ev_J$. Pour une vari\'et\'e convexe tous les arguments pr\'ec\'edents fonctionnent \emph{\`a $J$ fix\'e} en choisissant $\ux$ g\'en\'eriquement, y compris en pr\'esence d'une structure r\'eelle~\cite{We-Duke}.

\subsection{L'invariant de Welschinger en dimension $3$} \label{sec:dim3}

On suppose dans cette section $n=3$ et $c_1(X)d$ pair. Le nombre de points marqu\'es n\'ecessaires pour rigidifier l'espace de modules de courbes $J$-holomorphes homologues \`a $d$ avec des conditions d'incidence ponctuelles est $k=k_d=\frac 1 2 c_1(X)d$. \'Etant donn\'e $J\in\R\cJ_\omega$ g\'en\'erique, la d\'efinition des signes de Welschinger des courbes $C\in\cR^d(\ux,J)$ se fait en plusieurs \'etapes. 

\smallskip 

\emph{Choix d'une structure $\Pin_3^-$ sur $\R X$.} Choisissons une m\'etrique riemannienne sur $\R X$. Le rev\^etement universel (\`a deux feuillets) du groupe orthogonal $O_3(\R)$ admet deux structures de groupe diff\'erentes, not\'ees $\Pin_3^\pm$, pour lesquelles la projection est un morphisme de groupes~: pour $\Pin_3^-$ (resp. $\Pin_3^+$) le relev\'e d'une r\'eflection est d'ordre quatre (resp. d'ordre deux). L'obstruction \`a relever le $O_3(\R)$-fibr\'e principal des rep\`eres sur $\R X$ \`a un $\Pin_3^-$-fibr\'e principal est donn\'ee par la classe caract\'eristique $w_2(\R X)+w_1^2(\R X)\in H^2(\R X;\Z/2\Z)$~\cite[Lemme~1.3]{KT}, o\`u $w_1(\R X)$, $w_2(\R X)$ sont les classes de Stiefel-Whitney. La formule de Wu pour la classe de Stiefel-Whitney totale~\cite[\S11]{Milnor} implique que l'obstruction s'annule pour une vari\'et\'e compacte de dimension trois. Une structure $\Pin_3^-$ est un $\Pin_3^-$-fibr\'e principal qui rel\`eve le fibr\'e  des rep\`eres de $\R X$. L'espace des structures $\Pin_3^-$ est affine sur $H^1(\R X;\Z/2\Z)$.

\begin{center}
\emph{On choisit par la suite une structure $\Pin_3^-$ sur $\R X$, not\'ee $\pp$.} 
\end{center}

Les composantes orientables de $\R X$ admettent des structures $\Spin_3$ (on a $w_1=0$ et l'obstruction \`a l'existence d'une structure $\Spin_3$ est $w_2=w_1^2=0$). Le choix d'une orientation de $\R X$ permet de r\'eduire la structure $\pp$ \`a une structure $\Spin_3$.

\smallskip

\emph{D\'efinition du signe de Welschinger d'un op\'erateur de Cauchy-Riemann g\'en\'eralis\'e r\'eel~\cite{We-semipositive}.} Soient $J\in\R\cJ_\omega$, $u$ une courbe $J$-holomorphe simple immerg\'ee, $\uz$ une collection r\'eelle de $k_d$ points distincts dans $\P^1$. Supposons que $u$ a une partie r\'eelle non-vide. Soit $N_u$  le fibr\'e normal de $u$, regard\'e en tant que fibr\'e complexe. Soit $Op_\dbar(N_u)$ l'espace des op\'erateurs de Cauchy-Riemann sur $N_u$, $\R Op_\dbar(N_u)$ l'espace des op\'erateurs de Cauchy-Riemann r\'eels (invariants sous l'action de $\Z/2\Z$), $Op_{\dbar+R}(N_u)$ l'espace des op\'erateurs de Cauchy-Riemann g\'en\'eralis\'es sur $N_u$, et $\R Op_{\dbar+R}(N_u)$ l'espace des op\'erateurs de Cauchy-Riemann g\'en\'eralis\'es r\'eels. Nous consid\'erons des op\'erateurs de classe $C^{\ell-1}$. \`A titre d'exemple, $Op_{\dbar+R}(N_u)$ est un espace affine sur $C^{\ell-1}(\Lambda^{0,1}\P^1\otimes \mathrm{End}_\R(N_u))$, alors que $\R Op_{\dbar+R}(N_u)$ est un espace affine sur $C^{\ell-1}(\Lambda^{0,1}\P^1\otimes \mathrm{End}_\R(N_u))_{+1}$, le sous-espace propre correspondant \`a la valeur propre $+1$ pour l'action de $\Z/2\Z$. On note $D_\uz$ la restriction de $D$ \`a $W^{1,p}_{-\uz}(N_u)$. On d\'efinit le \emph{signe de Welschinger\footnote{Welschinger appelle ce signe \emph{\'etat spinoriel de $D$} et le note $sp(D)$~\cite{We-Duke,We-semipositive}.}} $\eps^\pp(D)\in\{\pm1\}$ d'un op\'erateur $D\in \R Op_{\dbar+R}(N_u)$ tel que $D_\uz$ soit inversible en trois \'etapes. 

\noindent (i)  On d\'efinit le signe de Welschinger $\eps^\pp(D)\in\{\pm1\}$ d'un op\'erateur $D\in \R Op_\dbar(N_u)$ tel que $D_\uz$ soit inversible de la mani\`ere suivante. L'op\'erateur $D$ d\'efinit une structure holomorphe sur $N_u$ pour laquelle on a une d\'ecomposition $N_u\simeq \cO(k_d-1)\oplus \cO(k_d-1)$ (Corollaire~\ref{coro:balanced}). Cette d\'ecomposition est compatible avec l'action de $\Z/2\Z$ et on en d\'eduit un scindement $\R N_u=L\oplus M$. Notons $C=\mathrm{im}(u)$, de sorte que $\R C$ d\'efinit un n\oe ud immerg\'e dans $\R X$. Ce n\oe ud est \'equip\'e d'un rep\`ere d'axes mobiles $T\R C\oplus L \oplus M$. Supposons que $L$ et $M$ sont orientables, i.e. $k_d-1$ est pair. Le rep\`ere d'axes mobiles peut \^etre alors enrichi en un rep\`ere mobile qui rel\`eve $\R C$ au fibr\'e des rep\`eres de $\R X$. On d\'efinit $\eps^\pp(D)=\pm1$ selon que ce lacet se rel\`eve au $\Pin_3^-$-fibr\'e principal $\pp$, ou pas. Dans le cas o\`u $L$ et $M$ ne sont pas orientables, on les tord par un demi-tour \`a droite le long du fibr\'e trivial $T\R C$, et on se ram\`ene \`a la situation pr\'ec\'edente~\cite[\S2.2]{We-Duke}~; 

\noindent (ii) On montre que l'espace $\R Op_{\dbar+R}(N_u)^{sing}\subset \R Op_{\dbar+R}(N_u)$ des op\'erateurs tels que $D_\uz$ n'est pas inversible est contenu dans une union d\'enombrable de sous-vari\'et\'es de codimension $\ge 1$~\cite[\S3.1]{We-semipositive}~; 

\noindent (iii) On d\'efinit le signe de Welschinger $\eps^\pp(D)\in\{\pm1\}$ d'un op\'erateur $D\in \R Op_{\dbar+R}(N_u)$ tel que $D_\uz$ soit inversible de la mani\`ere suivante. On choisit un chemin g\'en\'erique $\gamma$ dans $\R Op_{\dbar+R}(N_u)$ reliant $D$ \`a $D'\in \R Op_\dbar(N_u)$ avec $D'_\uz$ inversible, on note $n(\gamma)$ le nombre de fois que le chemin $\gamma$ intersecte le mur des op\'erateurs $D^{sing}$ tels que $D^{sing}_\uz$ a un conoyau de dimension $1$, et on pose 
$$
\eps^\pp(D)=(-1)^{n(\gamma)}\eps^\pp(D').
$$
Un argument d'intersection montre que la valeur $\eps^\pp(D)$ ne d\'epend pas du choix de $\gamma$. Par ailleurs, Welschinger prouve que le r\'esultat ne d\'epend pas du choix de $D'$ non-plus. Ceci revient \`a montrer que, lorsque $D'$ et $D''$ sont adjacents \`a un m\^eme mur d'op\'erateurs ayant un conoyau de dimension $1$, leurs signes diff\`erent~\cite[Proposition~3.2]{We-semipositive}. Le ph\'enom\`ene sous-jacent est le suivant~: au moment de la travers\'ee du mur le fibr\'e normal $N_u$ scinde comme $N_u=\cO(k_d)\oplus \cO(k_d-2)$. Ceci correspond au fait que, dans le cours de la d\'eformation, l'une des sections qui engendraient un sommand direct traverse la section nulle et acquiert un z\'ero suppl\'ementaire. Ceci implique que le lacet des rep\`eres d'un c\^ot\'e du mur se d\'eduit du lacet de l'autre c\^ot\'e en rajoutant un g\'en\'erateur de $\pi_1(SO(3))$. Un et un seul de ces deux lacets se rel\`eve donc au fibr\'e $\pp$, ce qui \'equivaut \`a dire que les signes de Welschinger de $D'$ et $D''$ sont oppos\'es.  

\smallskip 

\emph{D\'efinition du signe de Welschinger d'une courbe $C\in\cR^d(\ux,J)$.} Soit $\ux\in\R(X^k\setminus\mathrm{Diag})$ fix\'e et $J\in\R\cJ_\omega$ assez g\'en\'erique pour que les \'el\'ements de $\cR^d(\ux,J)$ soient des courbes plong\'ees et des valeurs r\'eguli\`eres de l'application d'\'evaluation $\R\ev_J:\R\cM^d_k(X,J)\to\R(X^k)$, en particulier des valeurs r\'eguli\`eres de l'application d'\'evaluation $\ev_J:\cM^d_k(X,J)\to X^k$. Chaque telle courbe d\'etermine un op\'erateur $D=\oD_u$ sur $N_u$ dont la restriction $D_\uz$ est un isomorphisme (Corollaire~\ref{coro:balanced}). On pose 
$$
\eps^\pp(C)=\eps^\pp(\oD_u).
$$

\smallskip 

\noindent {\it Esquisse de la d\'emonstration du th\'eor\`eme~\ref{thm:main} dans le cas $n=3$~\cite{We-Duke,We-semipositive}.} --- La condition $r\ge 1$ impose que la partie r\'eelle des courbes consid\'er\'ees est non-vide, de sorte que leurs signes de Welschinger sont d\'efinis. La d\'emonstration de l'invariance de $\chi^d_r(\ux,J)$ par rapport au choix de $\ux$ et $J$ suit exactement les m\^emes \'etapes que dans le cas $n=2$. Le cas des courbes cuspidales est remplac\'e par le cas des courbes telles que $D_\uz$ a un conoyau de dimension $1$. Ce cas est tautologique, au sens o\`u il a d\'ej\`a d\^u \^etre trait\'e pour montrer l'ind\'ependance de $\eps^\pp(D)$, $D\in\R Op_{\dbar+R}(N_u)$ par rapport au choix de $D'\in Op_\dbar(N_u)$, cf. ci-dessus. Le cas des courbes r\'eductibles est trait\'e par Welschinger en r\'eduisant le probl\`eme \`a la situation o\`u $J$ est int\'egrable au voisinage de la courbe, une situation similaire \`a celle rencontr\'ee dans le cas des vari\'et\'es convexes. Cette simplification est rendue possible par le fait que l'on peut choisir \`a volont\'e le point par lequel on traverse le mur des courbes r\'eductibles, ainsi que le segment $J_t$ avec lequel on traverse ce mur. Le probl\`eme \`a $\ux$ fix\'e et $J$ variable est alors \'equivalent \`a un probl\`eme \`a $J$ fix\'e et $\ux$ variable. Le fait que $J$ soit fix\'e et int\'egrable entra\^{\i}ne que la structure holomorphe sur le fibr\'e $N_u$ est constante. La preuve est finie en analysant le comportement d'une section holomorphe appropri\'ee lorsque $\ux$ franchit le mur~\cite[Proposition~3.7]{We-Duke}.  
\hfill{$\square$}

\smallskip

\noindent {\it Remarque (signes en dimension $n\ge 4$).} --- La d\'efinition des signes de Welschinger en dimension $n\ge 4$ suit les m\^emes \'etapes qu'en dimension trois~: (i) on impose des conditions topologiques sur $\R X$ qui garantissent l'existence d'une structure $\Pin_n^\pm$ ou $\Spin_n$ sur un certain fibr\'e d\'efini \`a partir de $T\R X$ et contenant ce dernier comme sous-fibr\'e~; (ii) on choisit une telle structure $\pp$ et on d\'efinit un signe pour les op\'erateurs de Cauchy-Riemann r\'eels sur $N_u$ tels que $D_\uz$ est inversible~; (iii) on en d\'eduit un signe pour les op\'erateurs de Cauchy-Riemann r\'eels g\'en\'eralis\'es sur $N_u$ tels que $D_\uz$ soit inversible~; (iv) \'etant donn\'ee $C\in\cR^d(\ux,J)$, on d\'efinit le signe $\eps^\pp(C)=\eps^\pp(\oD_u)$. 
\`A la diff\'erence du cas $n=3$, le nombre $\chi_r^{d,\pp}(\ux,J)$ n'est plus invariant lors de la travers\'ee de certains murs de courbes r\'eductibles.  
\emph{C'est une question ouverte importante que de comprendre quelle est la bonne d\'efinition d'un invariant de Welschinger en dimension $n\ge 4$.} En particulier, il me semble important de donner une d\'efinition  purement analytique des signes de Welschinger pour $n=3$.

\smallskip 

\noindent {\it Remarque (signes en dimension $2$).} \label{rmq:signes2} --- Les signes de Welschinger en dimension $n=2$ peuvent \^etre reformul\'es dans le langage de cette section. Supposons pour simplifier que $\R X$ est orientable et soit $\ss$ une structure $\Spin_2$ sur $T\R X$. Fixons une orientation sur $\R X$. \'Etant donn\'ee $C\in\cR^d(\ux,J)$, les trois informations suivantes sont \'equivalentes~: 
\begin{itemize}
\item la parit\'e de la masse de $C$~; 
\item la parit\'e du nombre de n\oe uds r\'eels non-solitaires de $C$~; 
\item la parit\'e du nombre de tours effectu\'es par le lacet orient\'e $\R C$, mesur\'ee par la structure $\ss$ selon la recette suivante~: soit $\tau$ le champ tangent orient\'e le long de $\R C$ et $v$ un champ le long de $\R C$ tel que $(\tau,v)$ soit une base positivement orient\'ee de $T\R X$ en chaque point de $\R C$. On dit que $\R C$ effectue un nombre pair, resp. impair de tours par rapport \`a $\ss$ si le lacet $(\tau,v)$ \`a valeurs dans le fibr\'e des rep\`eres de $T\R X$ se rel\`eve au fibr\'e $\ss$, resp. ne se rel\`eve pas.   
\end{itemize}

\section{Optimalit\'e, congruences} \label{sec:optimalite}

Dans cette section nous abordons deux questions essentielles li\'ees aux invariants $\chi^d_r$ du~\S\ref{sec:chird} en dimension quatre. 
La premi\`ere question est celle de l'optimalit\'e des bornes inf\'erieures fournies par $|\chi^d_r|$, cf. Corollaire~\ref{cor:borneinf}.

\begin{center}
\emph{Existe-t-il des structures presque complexes r\'eelles g\'en\'eriques telles que $|\chi^d_r|=R^d(\ux,J)$~?} 
\end{center}

Cela revient \`a demander s'il existe des structures presque complexes r\'eelles g\'en\'eriques pour lesquelles les \'el\'ements de $\cR^d(\ux,J)$ sont compt\'es avec le m\^eme signe. Un exemple de telle situation, qui appara\^{\i}tra plus bas, est celui o\`u tous les \'el\'ements de $\cR^d(\ux,J)$ ont un lieu r\'eel plong\'e. Les \'eventuels n\oe uds r\'eels d'une courbe $C$ sont alors solitaires et leur nombre, \'egal \`a la masse $m(C)$, est de m\^eme parit\'e que le nombre total $\delta(C)$ de points doubles de $C$, encore appel\'e \emph{le genre lisse} de $C$. La courbe $C$ \'etant immerg\'ee, on a $\delta(C)=\frac 1 2 (d^2-c_1(X)d+2)$ par la formule d'adjonction. Tous les \'el\'ements de $\cR^d(\ux,J)$ sont donc compt\'es dans cette situation avec le m\^eme signe $\eps=(-1)^{\frac 1 2 (d^2-c_1(X)d+2)}$.

\smallskip
 
La deuxi\`eme question est celle des congruences~: 

\begin{center}
\emph{Peut-on identifier dans certaines situations des (grands) diviseurs de $\chi^d_r$~?} 
\end{center}

Une r\'eponse affirmative a des cons\'equences imm\'ediates sur la borne inf\'erieure $|\chi^d_r|$~: si $\chi^d_r$ n'est pas nul et est divisible par un entier $m\ge 2$, alors $R^d(\ux,J)\ge m$ pour tous $\ux$ et $J$.
Mais l'enjeu ne s'arr\^ete pas l\`a~: les probl\`emes de congruence touchent historiquement au c\oe ur de la g\'eom\'etrie alg\'ebrique r\'eelle et sont la manifestation de ph\'enom\`enes profonds. \`A titre d'exemple, \'evoquons la congruence de Gudkov-Arnol'd-Rokhlin qui a marqu\'e l'irruption des m\'ethodes de topologie des vari\'et\'es de dimension quatre en g\'eom\'etrie alg\'ebrique r\'eelle dans les ann\'ees 1970~\cite[\S1]{DK}~: \emph{soit $C$ une courbe plane r\'eelle maximale de degr\'e $2k$, i.e. ayant un nombre maximal de composantes connexes (ovales). Soit $p$, resp. $n$, le nombre d'ovales contenus \`a l'int\'erieur d'un nombre pair, resp. impair d'autres ovales. Alors $p-n \equiv k^2 \ (mod \ 8)$.} 

Les r\'esultats de Welschinger que nous pr\'esentons dans cette section marquent l'irruption des m\'ethodes de la th\'eorie symplectique des champs en g\'eom\'etrie alg\'ebrique r\'eelle. 

\subsection{\'Enonc\'es}

\begin{theo} (optimalit\'e~\cite[Th\'eor\`eme~1.1]{We-optimalite}) \label{thm:optimalite} Supposons que $\R X$ contient une sph\`ere ou un plan projectif r\'eel $L$, et dans ce deuxi\`eme cas supposons que $(X,\omega,c_X)$ est symplectomorphe au plan projectif complexe \'eclat\'e en six points complexes conjugu\'es au maximum. Soit $0\le r\le 1$. Il existe des structures presque complexes r\'eguli\`eres telles que $|\chi^d_r(L)|=R^d(\ux,J)$. 
\end{theo} 

\smallskip 

\noindent {\it Remarque (signe).} --- Nous allons exhiber dans la preuve des $J$ r\'eguliers tels que les \'el\'ements de $\cR^d(\ux,J)$ ont tous un lieu r\'eel plong\'e. En vue de la discussion pr\'ec\'edente le signe de $\chi^d_r$ sera alors d\'etermin\'e par l'in\'egalit\'e $(-1)^{\frac 1 2 (d^2-c_1(X)d+2)}\chi^d_r\ge 0$.

\smallskip

\noindent {\it Remarque (maximalit\'e).} --- En \'ecrivant $k_d=c_1(X)d-1=r+2r_X$, la condition $0\le r\le 1$ peut \^etre formul\'ee de fa\c{c}on \'equivalente comme une condition de maximalit\'e pour le nombre $r_X$ de paires de points complexes conjugu\'es dans $\ux$. Il est utile de comparer cette situation \`a celle oppos\'ee, o\`u $r=c_1(X)d-1$ est maximal, qui est le cadre des r\'esultats de Itenberg, Kharlamov et Shustin mentionn\'es dans la section \S\ref{sec:developpements}. 

\smallskip 

\begin{coro} (\cite[Corollaire~1.3]{We-optimalite}) \label{coro:optimalite}
Soit $X$ le plan projectif complexe ou la quadrique ellipso\"{\i}de de dimension deux, et choisissons $0\le r\le 1$. On a l'\'egalit\'e $|\chi^d_r|=R^d(\ux,J)$ pour la structure complexe standard lorsque les points complexes conjugu\'es sont choisis tr\`es proches d'une conique imaginaire pure dans le premier cas et d'une section hyperplane r\'eelle disjointe de $L$ dans le second. 
\end{coro} 

\smallskip 

\noindent {\it Remarque.} --- Nous renvoyons au papier~\cite[\S1]{We-optimalite} pour deux autres r\'esultats d'optimalit\'e~: en dimension quatre, lorsque $r=1$ et le lieu r\'eel contient un tore, et en dimension six pour certaines vari\'et\'es convexes dont le lieu r\'eel contient une sph\`ere.  


%

\smallskip 

\begin{theo} (congruence~\cite[Th\'eor\`eme~2.1]{We-optimalite}) \label{thm:congruence}
Supposons que $\R X$ contient une sph\`ere ou un plan projectif r\'eel $L$, et dans ce deuxi\`eme cas supposons que $(X,\omega,c_X)$ est symplectomorphe au plan projectif complexe \'eclat\'e en six points complexes conjugu\'es au maximum. Si $L=S^2$ et $2r+1<k_d$, la puissance $2^{\frac 1 2 (k_d-2r-1)}=2^{\frac 1 2 (c_1(X)d - 2r-2)}$ divise $\chi^d_r(L)$. Si $L=\R P^2$ et $r+1<r_X$, la puissance $2^{r_X-r-1}$ divise $\chi^d_r(L)$.
\end{theo}

\smallskip

\noindent {\it Remarque.} --- Nous renvoyons le lecteur \`a~\cite[\S2]{We-optimalite} pour des variantes raffinn\'ees de ce th\'eor\`eme, ainsi que pour un autre r\'esultat de congruence concernant la quadrique ellipso\"{\i}de de dimension trois.

\subsection{Th\'eorie symplectique des champs} \label{sec:SFT}

L'ingr\'edient conceptuel nouveau qui intervient dans les preuves des th\'eor\`emes~\ref{thm:optimalite} et~\ref{thm:congruence} est la \emph{th\'eorie symplectique des champs} invent\'ee par Eliashberg, Givental et Hofer~\cite{EGH}. L'outil-cl\'e de celle-ci est le \emph{th\'eor\`eme de compacit\'e}~\cite{BEHWZ}. Soit $u_\nu$ une suite de courbes $J_\nu$-holomorphes \`a valeurs dans $X$, avec $J_\nu\in \cJ_\omega$. Le th\'eor\`eme de compacit\'e de Gromov~\cite{Gromov,MS04} d\'ecrit les d\'eg\'en\'erescences possibles dans la suite $u_\nu$ lorsque $J_\nu$ tend vers une limite $J_\infty\in\cJ_\omega$ avec $\nu\to\infty$. Le th\'eor\`eme de compacit\'e en th\'eorie symplectique des champs d\'ecrit les d\'eg\'en\'erescences possibles dans la suite $u_\nu$ lorsque la suite $J_\nu$ acquiert une singularit\'e d'un type tr\`es particulier lorsque $\nu\to\infty$, singularit\'e qui est concentr\'ee le long d'une hypersurface de type contact $\Sigma^{2n-1}\subset X^{2n}$. Cette condition signifie que $\omega$ admet au voisinage de $\Sigma$ une primitive $\lambda$ telle que $\alpha=\lambda|_\Sigma$ est une forme de contact, i.e. $\alpha\wedge d\alpha^{n-1}\neq 0$. En pr\'esence d'une structure r\'eelle, on demande que $\Sigma$ soit $c_X$-invariante et $c_X^*\alpha=-\alpha$. Un deuxi\`eme outil-cl\'e de la th\'eorie symplectique des champs est le th\'eor\`eme de recollement~\cite{EGH,HWZ99,B,HWZ09}. Sous des hypoth\`eses de transversalit\'e appropri\'ees celui-ci assure que toute configuration de courbes holomorphes qui est une limite possible est aussi une limite v\'eritable.

Nous expliquons dans cette section les ph\'enom\`enes qui se rattachent \`a la th\'eorie symplectique des champs dans notre contexte, \`a savoir celui d'une paire $(X,L)$ avec $L$ une composante connexe de $\R X$. Rappelons qu'un voisinage de $L$ dans $(X,\omega,c_X)$ est isomorphe \`a un voisinage de la section nulle dans $(T^*L,d\p\wedge d\q,c_L)$. Choisissons une m\'etrique riemannienne sur $L$ telle que $DT^*L=\{(\p,\q)\, : \, |\p|\le 1\}$ soit contenu dans ce voisinage et posons $\Sigma=ST^*L=\{(\p,\q)\, : \, |\p|=1\}$. La restriction $\alpha=\p d\q|_\Sigma$ est une forme de contact dont le champ de Reeb, d\'efini par $d\alpha(R,\cdot)=0$ et $\alpha(R)=1$, est le g\'en\'erateur du flot (co)g\'eod\'esique sur $ST^*L$. Le dual de $\p d\q$ par rapport \`a $d\p\wedge d\q$ est le champ de Liouville $\partial/\partial \p$ et son flot $\phi_t$ v\'erifie $\phi_t^*\p d\q=e^t\p d\q$. On obtient un symplectomorphisme $(x,t)\mapsto \phi_t(x)$ entre un cylindre $([-\eps,\eps]\times \Sigma,d(e^t\alpha))$, $\eps>0$ et un voisinage de $\Sigma$. Nous pouvons maintenant pr\'eciser les structures presque complexes singuli\`eres le long de $\Sigma$ qui sont admises par la th\'eorie symplectique des champs.
 
\begin{defi} \cite[\S2.1]{We-strings} Soit $J$ une structure presque complexe (r\'eelle) d\'efinie dans $X\setminus \Sigma$ et compatible avec $\omega$. On dit que \emph{$J$ est $\Sigma$-singuli\`ere} si 
\renewcommand{\theenumi}{\roman{enumi}}
\begin{enumerate}
\item $J$ pr\'eserve la distribution de contact $\xi=\ker\alpha$ pour chaque $\{t\}\times \Sigma$, $t\in[-\eps,\eps]\setminus \{0\}$ et sa restriction \`a $\xi$ ne d\'epend pas de $t\in[-\eps,\eps]\setminus \{0\}$~;
\item $J$ v\'erifie $J(\frac\partial {\partial t})=\beta'(t)R$ le long de $\{t\}\times \Sigma$, avec $\beta':[-\eps,\eps]\setminus\{0\}\to\R^*_+$ une fonction paire d'int\'egrale infinie.  
\end{enumerate}
\end{defi}

\begin{defi}  \cite[\S2.1]{We-strings} Soit $J$ une structure presque complexe (r\'eelle) d\'efinie dans $X$ et compatible avec $\omega$. On dit que \emph{$J$ a un cou cylindrique sur $\Sigma$} si 
\renewcommand{\theenumi}{\roman{enumi}}
\begin{enumerate}
\item $J$ pr\'eserve la distribution de contact $\xi=\ker\alpha$ pour chaque $\{t\}\times \Sigma$, $t\in[-\eps,\eps]$ et sa restriction \`a $\xi$ ne d\'epend pas de $t\in[-\eps,\eps]$~;
\item $J$ v\'erifie $J(\frac\partial {\partial t})=\beta'(t)R$ le long de $\{t\}\times \Sigma$, avec $\beta':[-\eps,\eps]\to\R^*_+$ une fonction paire.  
\end{enumerate}
La valeur de l'int\'egrale $\int_{-\eps}^\eps \beta'(t)dt$ est appel\'ee \emph{longueur du cou cylindrique}. 
\end{defi} 

\smallskip

\noindent {\it Remarque (structures presque complexes cylindriques).} --- Soit $J$ ayant un cou cylindrique de longueur $2A$ sur $\Sigma$ et notons $\beta:[-\eps,\eps]\to[-A,A]$ la primitive impaire de $\beta'$. Le pouss\'e en avant de $J$ par le diff\'eomorphisme $[-\eps,\eps]\times\Sigma\to[-A,A]\times\Sigma$, $(t,x)\mapsto(\beta(t),x)$ est une structure presque complexe $J'$ qui laisse invariante la distribution de contact, qui est invariante par translation en la variable $t$ et qui v\'erifie $J'(\frac\partial {\partial t})=R$. Ceci est par d\'efinition une structure presque complexe ``cylindrique'' au sens de la th\'eorie symplectique des champs~\cite[\S2]{BEHWZ}. Lorsque la structure presque complexe $J$ est $\Sigma$-singuli\`ere, la vari\'et\'e $X\setminus \Sigma$ est une vari\'et\'e \`a deux composantes $\mathrm{int}(DT^*L)$, resp. $X\setminus DT^*L$, chacune ayant un bout cylindrique semi-infini qui peut \^etre model\'e comme ci-dessus sur $\R_+\times\Sigma$, resp. $\R_-\times\Sigma$. Les courbes $J$-holomorphes \`a image dans $\mathrm{int}(DT^*L)$, resp. $X\setminus DT^*L$ seront identifi\'ees \`a des courbes $J'$-holomorphes \`a image dans $T^*L$, respectivement $X\setminus L$.

\smallskip

On se donne d\'esormais une suite $J_\nu$ de structures presque complexes ayant des cous cylindriques de longueur $2A_\nu\to\infty$ sur $\Sigma$, et on suppose que $J_\nu$ converge -- dans le sens \'evident -- vers une structure $\Sigma$-singuli\`ere $J_\infty$. Cette situation est appel\'ee de fa\c{c}on informelle \emph{\'etirement du cou}. Nous allons d\'ecrire maintenant les configurations de courbes qui sont limites de suites $u_\nu\in\cM^d_k(X,J_\nu)$, $d\in H_2(X;\Z)$, $k=k_d$. 

Notons $(W^+,J)$, resp. $(W^-,J)$ la vari\'et\'e $T^*L$, resp. $X\setminus L$, munie de la restriction de $J'_\infty$. Etant donn\'ee une courbe $J$-holomorphe $u:\dot S\to W^\pm$ d\'efinie sur une surface \'epoint\'ee, il existe une notion d'\emph{\'energie de Hofer}~(\cite[\S3.2]{Hofer}, \cite[\S5.3]{BEHWZ}) dont la finitude assure que $u$ est une application propre qui, en chaque pointe essentielle, est asymptote \`a un cylindre $Cyl_\gamma$ sur une orbite de Reeb p\'eriodique $\gamma$, de la forme $Cyl_\gamma(s,\theta)=(Ts+s_0,\gamma(T\theta))$ avec $(s,\theta)\in\R_\pm\times \R/2\pi\Z$~\cite{HWZ96,B}. Soulignons au passage cette vertu remarquable de l'\'energie de Hofer qui est celle de relier la dynamique hamiltonienne \`a la g\'eom\'etrie des courbes holomorphes. 

Dans le contexte qui nous int\'eresse, le th\'eor\`eme de compacit\'e en th\'eorie symplectique des champs prend la forme suivante. 

\begin{theo} \cite{BEHWZ,B} Soit $u_\nu\in\cM^d_k(X,J_\nu)$. Il existe une sous-suite, not\'ee $u_\nu$, qui converge au sens suivant vers une paire $u^\pm:\dot S^\pm\to W^\pm$ de courbes $J$-holomorphes d'\'energie de Hofer finie dont les pointes sont appari\'ees et qui ont les m\^emes orbites
asymptotes aux pointes correspondantes~: regardons $\dot S=\dot S^-\cup\dot S^+$ comme une courbe nodale en identifiant les pointes correspondantes de $\dot S^\pm$. Alors~:  
\begin{itemize}
\item il existe une suite d'applications $\phi_\nu:\P^1\to \dot S$ qui sont des diff\'eomorphismes en dehors d'une collection de cercles disjoints qui sont contract\'es sur les n\oe uds de $\dot S$. Ces cercles \'evitent les points marqu\'es sur $\P^1$~;  
\item $u_\nu\circ\phi_\nu^{-1}:\dot S^\pm\to X$ converge uniform\'ement sur tout compact vers $u^\pm$~;
\item les pointes de $u^\pm$ sont ``en phase"~: pour chaque n\oe ud $p$ de $\dot S$ consid\'erons une suite de segments $\gamma_\nu:]-\eps,\eps[\to \P^1$ qui intersectent $\phi_\nu^{-1}(p)$ transversalement en $s=0$ et tels que $\phi_\nu\circ \gamma_\nu=\gamma$. Alors $\lim_{s\to 0^+}\pi_\Sigma u^+(\gamma(s))=\lim_{s\to 0^-}\pi_\Sigma u^-(\gamma(s))$, o\`u $\pi_\Sigma u^\pm$ est la composante de $u^\pm$ sur $\Sigma$ dans la carte $\R_\pm\times \Sigma$. 
\end{itemize}
\end{theo}

\smallskip 

\noindent {\it Remarque (courbes \`a deux \'etages).} --- On appellera \emph{courbe \`a deux \'etages (dans $T^*L$ et $X\setminus L$)} un couple $(u^+,u^-)$ comme ci-dessus. Le th\'eor\`eme de compacit\'e pr\'edit l'existence d'\'etages interm\'ediaires qui sont des courbes holomorphes \`a valeurs dans la symplectisation $(\R\times \Sigma,d(e^t\alpha),J')$. Dans notre situation les courbes $u_\nu$ sont rigides et de tels \'etages interm\'ediaires n'apparaissent pas. Le th\'eor\`eme de compacit\'e pr\'edit aussi des courbes $u^-$ qui peuvent \^etre r\'eductibles, constitu\'ees d'une composante principale \'epoint\'ee \`a laquelle sont rattach\'ees des sph\`eres holomorphes. Ces courbes vivent en codimension $\ge 2$ dans l'espace de modules, de sorte qu'elles n'apparaissent pas dans notre situation. 

Chaque composante de $(u^+,u^-)$ est rigide, au sens o\`u elle appartient \`a un espace de modules de courbes \'epoint\'ees soumises \`a des conditions homologiques, asymptotiques et d'incidence qui est de dimension $0$. Chaque composante de $(u^+,u^-)$ est par ailleurs immerg\'ee. Le genre de $\dot S$ \'etant nul, deux composantes distinctes de $(u^+,u^-)$ ont au plus une asymptote commune. 

\smallskip 

\noindent {\it Remarque (courbes r\'eelles)}.Ê--- Lorsque les structures presque complexes et les courbes holomorphes $u_\nu$ sont r\'eelles, la courbe limite $(u^+,u^-)$ l'est aussi. La structure r\'eelle sur $T^*L$, resp. $X\setminus L$ est donn\'ee par $c_L$, resp. la restriction de $c_X$. Puisque le lieu r\'eel de $\P^1$ est connexe et disconnecte $\P^1$, la courbe $(u^+,u^-)$ a exactement une composante qui est $c_L$-invariante, alors que toutes les autres composantes viennent en paires qui sont conjugu\'ees par $c_L$ ou $c_X$. 

\smallskip 

On suppose maintenant $L=S^n$, resp. $L=\R P^n$, et on la munit  d'une m\'etrique riemannienne \`a courbure constante \'egale \`a $1$. Les orbites de Reeb ferm\'ees sur $\Sigma=ST^*L$, consid\'er\'ees modulo param\'etrisation, sont group\'ees en familles non-d\'eg\'en\'er\'ees au sens de Morse-Bott, de dimension \'egale \`a $2(n-1)$. Les p\'eriodes des orbites sont \'egales \`a $2k\pi$, $k\in\N^*$, resp. $k\pi$, $k\in\N^*$. L'entier $k$ est la \emph{multiplicit\'e} de l'orbite.

Fixons $p>2$. Soit $\dot S^\pm$ une surface \'epoint\'ee de genre $0$ et $\{F_1,\dots,F_v\}$ une collection de familles d'orbites de Reeb p\'eriodiques, index\'ee par les $v$ pointes de $\dot S^\pm$. \emph{Pour simplifier les notations, on suppose dans ce paragraphe que $\dot S^\pm$ est connexe.} Fixons $\delta>0$ et consid\'erons la vari\'et\'e de Banach $\cB$ des applications $u^\pm:\dot S^\pm\to W^\pm$ qui sont de classe $W^{1,p}_{\mathrm{loc}}$ et telles qu'au voisinage de chaque pointe $v_i$ on ait $u^\pm-Cyl_{\gamma_i}\in W^{1,p}(e^{\delta |s|/p}dsd\theta)$, avec $\gamma_i\in F_i$, $(s,\theta)\in\R_\pm\times \R/2\pi\Z$ et $Cyl_{\gamma_i}(s,\theta)=(Ts+s_0,\gamma_i(T\theta))$ le cylindre trivial au-dessus de $\gamma_i$. On dit encore que $u^\pm$ est de classe $W^{1,p,\delta}$.
Soit $\cE\to \cB$ le fibr\'e de Banach de fibre $\cE_u=L^{p,\delta}(\Lambda^{0,1}\dot S^\pm,(u^\pm)^*TW^\pm)$. De mani\`ere analogue au \S\ref{sec:courbes}, les courbes holomorphes $u^\pm:\dot S^\pm\to W^\pm$ qui sont asymptotes aux pointes $v_i$ \`a des orbites ferm\'ees appartenant aux $F_i$ sont les z\'eros de la section $\dbar_J:\cB\to \cE$ pour $\delta>0$ assez petit. Le fait d'utiliser des espaces de Sobolev \`a poids exponentiels est n\'ecessaire pour que l'op\'erateur lin\'earis\'e soit de Fredhlom, en raison de la pr\'esence de d\'eg\'en\'erescences le long des espaces tangents aux $F_i$, le long du champ de Reeb, et le long de la coordonn\'ee verticale $\partial/\partial t$ dans la symplectisation $\R_\pm\times \Sigma$. En lin\'earisant le probl\`eme on obtient un op\'erateur de Cauchy-Riemann g\'en\'eralis\'e 
$$
D_{u^\pm}:V\times W^{1,p,\delta}((u^\pm)^*TW^\pm)\to L^{p,\delta}(\Lambda^{0,1}\dot S^\pm,(u^\pm)^*TW^\pm),
$$
o\`u $V$ est un espace de dimension $N=\sum_{i=1}^v (\dim\, F_i+2)$ engendr\'e par des sections support\'ees le long de directions de d\'eg\'en\'erescence ind\'ependantes pour chaque pointe. Soit $\chi=2-v$ la caract\'eristique d'Euler de $\dot S^\pm$. L'op\'erateur $D_{u^\pm}$ est de Fredholm et son indice est donn\'e par~\cite[\S5]{B} 
\begin{eqnarray*}
\ind\, D_{u^\pm} & = & n\chi + \mu^{tot} + \frac N 2 \\
& = &  2n + \mu^{tot}.
\end{eqnarray*}
La deuxi\`eme \'egalit\'e d\'ecoule de ce que $\dim\,F_i=2(n-1)$, de sorte que $N=2nv$. Ici $\mu^{tot}$ d\'esigne \emph{l'indice de Maslov total} de $u^\pm$, qui est le double de l'obstruction \`a \'etendre \`a $\dot S^\pm$ la trivialisation de $(u^\pm)^*TW^\pm$ donn\'ee par le flot de Reeb lin\'earis\'e au voisinage des pointes. 
Lorsque la transversalit\'e est r\'ealis\'ee (par exemple lorsque $u^\pm$ sont des courbes simples), la dimension de l'espace de modules de courbes $\cM_{u^\pm}$ dans lequel vit $u^\pm$ est 
$$
\dim\,\cM_{u^\pm}=\ind\,D_{u^\pm} - 6 + 2v. 
$$
\begin{prop} \cite[Proposition~1.13]{We-optimalite}, \cite[Th\'eor\`eme~3.1]{Viterbo} \label{prop:mu}
Soit $L=S^n$ ou $\R P^n$. Soit $u^+:\dot S^+\to W^+=T^*L$ une courbe d'\'energie de Hofer finie et genre $0$. Soit $k$ la multiplicit\'e totale de ses orbites de Reeb asymptotes. L'indice de Maslov $\mu^{tot}$ de $u^+$ est $2k(n-1)$ lorsque $L=S^n$, respectivement $k(n-1)$ lorsque $L=\R P^n$.    
\end{prop}

Cet \'enonc\'e doit \^etre lu comme affirmant l'\'egalit\'e entre l'indice de Maslov total $\mu^{tot}$ et l'indice de Morse total de la collection des g\'eod\'esiques ferm\'ees qui correspondent aux orbites de Reeb asymptotes.

\subsection{D\'emonstrations}

Welschinger~\cite[\S1.1.2]{We-optimalite} encode le quotient par $\Z/2\Z$ d'une courbe r\'eelle \`a deux \'etages $u=(u^+,u^-)$ par un arbre enracin\'e dont les ar\^etes sont d\'ecor\'ees d'entiers positifs. 
La racine $s_0$ repr\'esente le quotient de l'unique composante $c_L$-invariante, qui est un disque \'epoint\'e \`a bord sur $L$. Les autres sommets repr\'esentent le quotient d'une paire de composantes complexes conjugu\'ees. Chaque ar\^ete adjacente \`a un sommet repr\'esente une paire d'asymptotes conjugu\'ees de la (paire de) composante(s) correspondante(s) et l'entier positif qu'elle porte est la multiplicit\'e de ces orbites asymptotes. 
De cette mani\`ere, les composantes \`a valeurs dans $T^*L$, resp. $X\setminus L$ sont repr\'esent\'ees par des sommets \`a distance paire, resp. impaire de $s_0$. 

\begin{figure}[h]
\label{fig:arbre}
\centering
\includegraphics[scale=0.8]{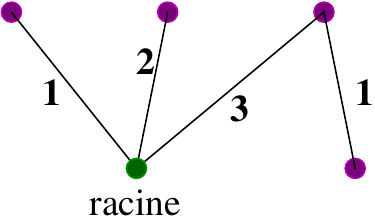}
\caption{Un exemple de courbe limite \`a $9$ composantes.}
\end{figure}

Les courbes limite $u^\pm$ sont immerg\'ees. Welschinger utilise \emph{l'indice de Maslov} $\mu$, d\'efini comme \'etant le double de l'obstruction \`a \'etendre \`a $\dot S^\pm$ la trivialisation du \emph{fibr\'e normal} $N_{u^\pm}$ donn\'ee par le flot de Reeb lin\'earis\'e au voisinage des pointes. En notant $\chi=2-v$ la caract\'eristique d'Euler d'une composante de $\dot S^\pm$ on obtient $\mu^{tot}=\mu+2\chi=\mu+4-2v$ et $\dim\cM_{u^\pm}=2n+\mu-2$. Lorsque $\dim\,X=4$ on a en particulier $\dim\,\cM_{u^\pm}=\mu+2$.

Cette derni\`ere formule de dimension est valable aussi lorsque l'on travaille avec des courbes r\'eelles, sauf pour la courbe correspondant au sommet $s_0$ de l'arbre, pour laquelle on a $\dim\cM_{s_0}=\frac 1 2 (\mu+2)=\frac 1 2 \mu +1$ puisqu'elle est $c_L$-invariante.

\smallskip

\noindent {\it D\'emonstration du th\'eor\`eme~\ref{thm:optimalite}.}
Nous pr\'esentons la preuve dans le cas $L=S^2$, le cas $L=\R P^2$ \'etant analogue. Remarquons que l'orientabilit\'e de $L$ implique l'imparit\'e de $r$~: soit $C$ un surface r\'eelle orientable immerg\'ee qui repr\'esente $d$ telle que $\R C\subset \R L$. Alors $r$ impair \'equivaut \`a $c_1(X)d$ pair, ou encore \`a ce que le fibr\'e normal de $C$ soit de degr\'e pair, ce qui d\'ecoule de l'orientabilit\'e du fibr\'e normal \`a $\R C$ dans $\R L$. 

On munit $L$ d'une m\'etrique \`a courbure constante et on \'etire le cou de la structure presque complexe au voisinage de $\Sigma=ST^*L$. L'id\'ee est de montrer que la courbe limite a un lieu r\'eel plong\'e. Ceci entra\^{\i}ne que, pour $\nu$ assez grand, les lieux r\'eels de $u_\nu$ sont plong\'es aussi, ce qui permet de conclure par l'argument pr\'esent\'e en d\'ebut de section. 

On note $A$ l'arbre qui encode la courbe limite et $S_1$, resp. $S_2$ l'ensemble des sommets \`a distance impaire, resp. paire de $s_0$. Pour chaque sommet $s$ on note $v_s$ sa valence et $k_s$ la somme des multiplicit\'es des ar\^etes adjacentes. On note $\mu_s$, resp. $\mu_s^{tot}$ les indices de Maslov de la courbe $C_s$ associ\'ee \`a un sommet $s$ et $\chi_s$ sa caract\'eristique d'Euler. On a en particulier $\chi_s=2-v_s$ pour $s\neq s_0$ et $\chi_s=1-v_s$ pour $s=s_0$. Soit $v$ le nombre total d'ar\^etes et $k$ leur multiplicit\'e totale. 

Pour les courbes de l'\'etage $T^*L$ on a $\mu_s=\mu_s^{tot}-2\chi_s$. Par la Proposition~\ref{prop:mu} on obtient 
$$
\sum_{s\in S_2} \mu_s=2k+2v-4\#S_2+2.
$$

Regardons maintenant les courbes de l'\'etage $X\setminus L$ et supposons pour commencer qu'elles sont simples. La g\'en\'ericit\'e de la structure presque complexe impose que tous les espaces de modules concern\'es sont de dimension positive, c'est-\`a-dire $\mu_s+2\ge 0$, respectivement $\mu_s+2\ge 2f_s$ si $C_s$ contient $f_s$ points de notre collection. On obtient la minoration 
$$
\sum_{s\in S_1} \mu_s \ge -2\#S_1+2r_X. 
$$

Puisque $A$ est un arbre on a $v=\#S_1+\#S_2-1$, de sorte que l'indice de Maslov $\mu=\sum_s\mu_s$ satisfait 
$$
\mu\ge 2k-2\#S_2+2r_X\ge 2r_X. 
$$
Par ailleurs $\mu$ est major\'e par $c_1(X)d-2$, le degr\'e du fibr\'e normal d'une courbe rationnelle immerg\'ee homologue \`a $d$. Puisque $r=1$ on a $c_1(X)d-2=2r_X$ et toutes les in\'egalit\'es pr\'ec\'edentes doivent \^etre des \'egalit\'es. En particulier $k=\#S_2$ , toutes les orbites de Reeb qui interviennent sont simples et tous les sommets de $S_2$ sont des feuilles, \emph{y compris $s_0$}. La courbe r\'eelle cod\'ee par $s_0$ est donc un cylindre r\'eel ayant comme asymptotes deux orbites de Reeb simples (conjugu\'ees). Welschinger d\'emontre par un raisonnement similaire \`a  celui qui prouve la formule d'adjonction qu'un tel cylindre est n\'ecessairement plong\'e~\cite[Lemme~1.14]{We-optimalite}, ce qui est la conclusion d\'esir\'ee. 

Le cas o\`u les courbes $C_s$ sont multiplement rev\^etues est trait\'e en utilisant les faits suivants~: (i) une courbe d'\'energie de Hofer finie factorise toujours \`a travers une courbe simple~\cite[Appendice]{HWZ95}~; (ii) l'indice de Maslov $\mu^\ell$ d'un rev\^etement de degr\'e $\ell$ d'une courbe simple d'indice $\mu$ vaut $\mu^\ell=\ell\mu+2\rho$, o\`u $\rho$ est l'indice de ramification. 
Il en d\'ecoule que $\mu^\ell$ peut \^etre plus petit que $\mu$ uniquement lorsque $\mu$ est n\'egatif, donc \'egal \`a $-2$. Cela ne concerne en particulier pas les courbes de l'\'etage $X\setminus L$ soumises \`a des conditions d'incidence, not\'ees $C_{s_1},\dots,C_{s_j}$, qui v\'erifient par cons\'equent $\sum_{i=1}^j\mu_{s_i}\ge -2j+2r_X$. Nous allons montrer la minoration $\sum_{s\notin\{s_1,\dots,s_j\}}\mu_s\ge 2j$, ce qui permettra alors de conclure comme pr\'ec\'edemment. 

Nous allons estimer la contribution \`a l'indice de Maslov total pour chaque composante connexe de $A\setminus\{s_1,\dots,s_j\}$.
Soit $A'$ une telle composante connexe et notons $S'_1$, resp. $S'_2$ l'ensemble de ses sommets qui, dans $A$, sont \`a distance impaire, resp. paire de $s_0$. On note $v'$, resp. $k'$ la valence totale, resp. la multiplicit\'e totale \emph{dans $A$} des sommets de $A'$. Comme pr\'ec\'edemment on obtient 
$$
\sum_{s\in S'_2} \mu_s=2k'+2v'-4\#S'_2+2\delta,
$$
o\`u $\delta$ vaut $1$ si $A'$ contient $s_0$ et $0$ sinon. Pour estimer $\sum_{s\in S'_1} \mu_s$ nous introduisons les notations suivantes concernant un sommet $s\in S'_1$~: on suppose que $C_s$ est un rev\^etement de degr\'e $\ell_s$ d'une courbe simple $\underline C_s$, avec indice de ramification $\rho_s$, on note $v_s,\underline v_s$ leurs nombres respectifs de pointes et $\chi_s=2-v_s$, $\underline \chi_s=2-\underline v_s$ leurs caract\'eristiques d'Euler. La formule de Riemann-Hurwitz assure l'\'egalit\'e $\ell_s\underline \chi_s = \chi_s + \rho_s$. Soit $k_s$ la multiplicit\'e totale des ar\^etes adjacentes \`a $s$. On obtient
\begin{eqnarray*}
\sum_{s\in S'_1} \mu_s & = & \sum_{s\in S'_1} (\ell_s\underline \mu_s + 2 \rho_s) \\
& \ge & -2\sum_{s\in S'_1} \ell_s + 2\sum_{s\in S'_1} (\ell_s\underline \chi_s -\chi_s) \\
& = & 2\sum_{s\in S'_1} (\ell_s-\ell_s\underline v_s + v_s) - 4\#S'_1 \\
& \ge & 2\sum_{s\in S'_1} (\ell_s-k_s+v_s) - 4\#S'_1. 
\end{eqnarray*}
Soient $v'_{int}$ le nombre total d'ar\^etes de $A'$ et $k'_{int}$ leur multiplicit\'e totale. Puisqu'il n'y a pas d'ar\^ete qui relie l'un des sommets $s_1,\dots,s_j$ \`a un sommet dans $S'_1$, on obtient que $\sum_{s\in S'_1}v_s=v'_{int}$ et $\sum_{s\in S'_1} k_s=k'_{int}$. Par ailleurs $v'_{int}=\#S'_1+\#S'_2-1$ et l'on obtient  
$$
\sum_{s\in S'_1\cup S'_2} \mu_s \ge 2(k'-k'_{int}) + 2(v'-v'_{int}) + 2\sum_{s\in S'_1} \ell_s - 4 + 2\delta.
$$
Dans le membre de droite $k'-k'_{int}\ge 1$, $v'-v'_{int}\ge 1$ et $\ell_s\ge 1$.
Chacun des sommets $s_1,\dots ,s_j$ \'etant reli\'e \`a au moins une composante $A'$ comme ci-dessus, pour laquelle il contribue de $1$ dans $k'-k'_{int}$ \emph{et} dans $v'-v'_{int}$, on obtient en sommant sur toutes les composantes $A'$ de $A\setminus\{s_1,\dots,s_j\}$ la minoration d\'esir\'ee
$$
\sum_{s\notin\{s_1,\dots,s_j\}}\mu_s\ge 2j.
$$ 
\hfill{$\square$}

\noindent {\it Remarque.} --- Il s'ensuit de la d\'emonstration que l'\'enonc\'e du Th\'eor\`eme~\ref{thm:optimalite}  peut \^etre pr\'ecis\'e~: les bornes inf\'erieures sont atteintes pour toute structure presque complexe g\'en\'erique ayant un cou suffisamment long au voisinage de $L$.

\smallskip

\noindent {\it D\'emonstration du Corollaire~\ref{coro:optimalite}.} --- Le plan projectif et la quadrique sont des surfaces convexes, de sorte que la structure complexe standard est g\'en\'erique. Dans les deux cas, la structure standard a un cou de longueur infinie au voisinage de $L$~: il s'agit du compl\'ementaire de la conique imaginaire pure, respectivement du compl\'ementaire de la section hyperplane. On conclut par la remarque pr\'ec\'edente.
\hfill{$\square$}

\smallskip

\noindent {\it Esquisse de la d\'emonstration du Th\'eor\`eme~\ref{thm:congruence}.} --- Il s'agit de d\'ecrire avec plus de d\'etail les arbres qui sont susceptibles d'encoder une courbe \`a deux \'etages qui est limite d'une suite $u_\nu\in\cM^d_k(X,J_\nu)$, $\nu\to\infty$ apr\`es \'etirement du cou. Toutes les composantes de la courbe limite sont rigidifi\'ees par leurs conditions d'incidence et leurs conditions asymptotiques. Dans le cas $L=S^2$, Welschinger montre en utilisant des estim\'ees sur l'indice de Maslov semblables \`a celles de la preuve du th\'eor\`eme~\ref{thm:optimalite} que toutes les composantes de l'\'etage $X\setminus L$ sont connect\'ees \`a la racine $C_{s_0}$, et qu'il y a au moins $\frac 1 2 (k_d-2r-1)$ paires complexes conjugu\'ees de telles composantes. Chaque paire est rigidifi\'ee en prescrivant une paire d'orbites asymptotes communes avec $C_{s_0}$. Puisqu'il y a deux mani\`eres d'apparier deux paires d'orbites et donc de recoller une telle composante de l'\'etage $X\setminus L$ \`a $C_{s_0}$, il en r\'esulte que la puissance $2^{\frac 1 2 (k_d-2r-1)}$ divise $\chi^d_r$. Lorsque $L=\R P^2$ le raisonnement est un peu plus d\'elicat parce-que les composantes de l'\'etage $X\setminus L$ ne sont pas n\'ecessairement connect\'ees \`a la racine $C_{s_0}$.

\subsection{Ouverture~: invariants relatifs} \label{sec:invariants-relatifs}

Revenons \`a la courbe \`a deux \'etages obtenue apr\`es \'etirement du cou au voisinage de la lagrangienne $L$, cod\'ee par un arbre enracin\'e $A$. Nous pouvons r\'esumer la d\'emarche suivie jusqu'ici de la mani\`ere suivante~: 
\begin{enumerate}
\item une description grossi\`ere de l'arbre $A$, n'utilisant essentiellement que la structure \`a deux \'etages de la courbe et des estim\'ees sur l'indice de Maslov, a permis d'obtenir le th\'eor\`eme d'optimalit\'e~\ref{thm:optimalite}~; 
\item une description plus fine de l'arbre $A$, avec des informations sur la distance des composantes de l'\'etage $X\setminus L$ par rapport \`a la racine $s_0$, permet d'obtenir le th\'eor\`eme de congruence~\ref{thm:congruence}~; 
\item une description \emph{exhaustive} du type combinatoire de $A$ permet d'exprimer $\chi^d_r$ comme un produit de convolution d'invariants d\'efinis dans $T^*L$ et d'invariants d\'efinis dans $X\setminus L$. La convolution est entendue ici comme une somme discr\`ete sur tous les types combinatoires possibles de courbes \`a deux \'etages. Avec une notation vague on peut \'ecrire  
$$
\chi^d_r=\chi_{T^*L} * \chi_{X\setminus L}.
$$
Ceci est le reflet alg\'ebrique du ``cassage" de la vari\'et\'e $X$ en deux morceaux $T^*L$ et $X\setminus L$ par \'etirement du cou.
\end{enumerate}

Welschinger a rendu  rigoureux ce dernier point dans~\cite[Th\'eor\`emes~3.10 et~3.16]{We-optimalite} lorsque $X$ est le plan projectif, respectivement la quadrique ellipso\"{\i}de de dimension $2$.  
\renewcommand{\theenumi}{\roman{enumi}}
\begin{enumerate}
\item Le cas $X=\P^2$. On travaille avec la structure complexe standard, auquel cas $T^*\R P^2$ devient biholomorphe \`a $\P^2\setminus Q$, o\`u $Q$ est la conique imaginaire pure, et $\P^2\setminus \R P^2$ devient biholomorphe \`a l'espace total du fibr\'e holomorphe de degr\'e quatre sur $Q$. Les courbes situ\'ees dans l'\'etage $T^*L$, soumises \`a des conditions d'incidence totalement r\'eelles et \`a des conditions asymptotiques, peuvent \^etre interpr\'et\'ees comme des invariants relatifs au diviseur r\'eel $Q$, la multiplicit\'e d'une orbite asymptote encodant l'ordre de tangence \`a $Q$. Ce type d'invariant a \'et\'e d\'efini par Welschinger dans~\cite{We07a}. Le comptage des courbes situ\'ees dans l'\'etage $X\setminus L$, ayant des conditions d'incidence complexes conjugu\'ees, peut \^etre interpr\'et\'e comme un invariant de Gromov-Witten relatif dans la surface rationnelle r\'eelle regl\'ee de degr\'e quatre, ayant des conditions de tangence prescrites avec la section exceptionnelle.
\item Le cas $X=\P^1\times \P^1$, la quadrique ellipso\"{\i}de de dimension deux.  On a $\R X=S^2$. On travaille avec la structure complexe standard, auquel cas $T^*S^2$ devient biholomorphe \`a $(\P^1\times \P^1)\setminus Q$, o\`u $Q$ est une section hyperplane r\'eelle de $X$ disjointe de $\R X = S^2$, et $(\P^1\times \P^1)\setminus S^2$ devient biholomorphe \`a l'espace total du fibr\'e holomorphe de degr\'e deux sur $Q$.  Les courbes situ\'ees dans l'\'etage $T^*L$, soumises \`a des conditions d'incidence totalement r\'eelles et \`a des conditions asymptotiques, peuvent \^etre interpr\'et\'ees comme des invariants relatifs au diviseur r\'eel $Q$~\cite{We07a}, la multiplicit\'e d'une orbite asymptote encodant l'ordre de tangence \`a $Q$. Le comptage des courbes situ\'ees dans l'\'etage $X\setminus L$, ayant des conditions d'incidence complexes conjugu\'ees, peut \^etre interpr\'et\'e comme un invariant de Gromov-Witten relatif dans la surface rationnelle r\'eelle regl\'ee de degr\'e deux, ayant des conditions de tangence prescrites avec la section exceptionnelle.
\end{enumerate}
En vue de la discussion pr\'ec\'edente, lorsque $X=\P^2$ ou $X=\P^1\times \P^1$ 
l'\'equation de convolution ci-dessus prend la forme
$$
\chi^d_r=\chi^{rel}_{T^*L} * GW^{rel}_{\overline{X\setminus L}}
$$
et exprime $\chi^d_r$ comme un produit de convolution entre un invariant de Welschinger relatif et un invariant de Gromov-Witten relatif dans une compactification appropri\'ee de $X\setminus L$. \`A nouveau, la convolution est entendue comme une somme discr\`ete sur tous les types combinatoires possibles de courbes \`a deux \'etages. La formule est de m\^eme nature que la formule de Ionel et Parker exprimant les invariants de Gromov-Witten d'une somme connexe symplectique le long d'un diviseur comme un produit de convolution d'invariants de Gromov-Witten relatifs au diviseur~\cite{IP}. 

\smallskip

\noindent {\it Remarque (autres invariants relatifs).}Ê--- Welschinger d\'efinit dans~\cite{We06} des invariants relatifs r\'eels de vari\'et\'es de dimension quatre en imposant des conditions de tangence au lieu r\'eel de diviseurs particuliers. Ceci lui a permis en particulier de montrer que le nombre de coniques r\'eelles tangentes \`a cinq coniques r\'eelles g\'en\'eriques de $\P^2$ est toujours minor\'e par $32$. De Joncqui\`eres avait \'etabli en 1859 que le nombre de coniques complexes vaut 3264, alors que Ronga, Tognoli et Vust avaient montr\'e en 1997 qu'il existe des configurations r\'eelles pour lesquelles toutes les solutions sont r\'eelles~\cite{Ronga}.

\section{Autres d\'eveloppements} \label{sec:developpements}

\subsection{Sym\'etrie miroir} \label{sec:symetrie-miroir}
 
Solomon~\cite[Th\'eor\`eme~1.3]{S06}{\,Ê} a d\'efini en dimensions $2$ et $3$, et sous l'hypoth\`ese que $\R X$ est orientable dans ce deuxi\`eme cas, des invariants \'enum\'eratifs r\'eels qui g\'en\'eralisent les invariants de Welschinger. Les invariants de Solomon comptent des courbes $J$-holomorphes soumises \`a des conditions d'incidence ponctuelles, en genre arbitraire et ayant un nombre arbitraire de composantes de bord, contraintes \`a avoir une image dans $\R X$. Dans l'approche de Solomon la structure conforme \`a la source est \emph{fix\'ee}. La m\'ethode de construction de ces invariants est proche de celle des invariants de Gromov-Witten et fournit en particulier une interpr\'etation des invariants de Welschinger $\chi^{d,\pp}_r(L)$ en termes d'int\'egrales de formes diff\'erentielles sur un espace de modules de disques $J$-holomorphes avec condition au bord lagrangienne~\cite[Th\'eor\`eme~1.8]{S06}. Cet espace de modules de disques est un rev\^etement double de l'espace de modules de courbes rationnelles r\'eelles. 

Une interpr\'etation des invariants de Welschinger dans le m\^eme esprit a \'et\'e donn\'ee par Cho~\cite{Cho} sous l'hypoth\`ese que $\R X$ est orientable. Dans la m\^eme direction, mentionnons~\cite{FOOO-paper} et l'article de Fukaya~\cite{F}, ainsi que l'ouvrage fondateur~\cite{FOOO-book}.

Tous ces travaux constituent autant d'approches au probl\`eme de d\'efinir des invariant de Gromov-Witten ``ouverts'', i.e. une th\'eorie d'intersection coh\'erente sur l'espace de modules d'applications stables avec conditions au bord lagrangienne. Il semble que les vari\'et\'es symplectiques r\'eelles fournissent un cadre appropri\'e pour ce probl\`eme fondamental, la structure r\'eelle assurant des annulations miraculeuses pour des termes de bord d'int\'egrales d\'efinies sur l'espace de modules. L'un des probl\`emes centraux du domaine est de d\'efinir des invariants de Welschinger en dimension $2n\ge 8$. 


Solomon a d\'ej\`a expos\'e des r\'esultats concernant un analogue de l'\'equation WDVV~\cite[\S11.2]{MS04} pour l'espace de modules de disques stables \`a bord lagrangien (cf. \cite{ABLdM}). Ceci sugg\`ere l'existence d'une version \emph{r\'eelle} de la conjecture de sym\'etrie miroir~\cite{K94}. 
Pandharipande, Solomon et Walcher \cite{PSW} calculent des invariants \'enum\'eratifs pour la quintinque r\'eelle de $\mathbb{P}^4$ en utilisant des formes d\'ej\`a d\'emontr\'ees de sym\'etrie miroir. (La quintique sort du cadre des vari\'et\'es semi-positives que nous avons adopt\'e ici.)

\noindent {\it Remarque (Invariants de Gromov-Witten ``ouverts'' en dimension quatre).}Ê--- Alors que j'\'etais en train de mettre la derni\`ere main \`a cet article avant publication, Welschinger vient de d\'efinir des invariants de Gromov-Witten ouverts \'enum\'eratifs en dimension quatre~\cite{We11}. Il s'agit d'un comptage avec signe de disques \`a bord sur une sous-vari\'et\'e lagrangienne orientable, sous la seule condition que le bord soit homologiquement trivial et \emph{sans utiliser de structure r\'eelle}. Ces nouveaux invariants g\'en\'eralisent ceux du \S\ref{sec:chird}.

\subsection{G\'eom\'etrie tropicale}

La g\'eom\'etrie tropicale peut \^etre d\'ecrite comme \'etant la g\'eom\'etrie alg\'ebrique sur le semi-anneau tropical $\R_{trop}=(\R,\max,+)$. Les op\'erations $\max$ et $+$ peuvent \^etre vues comme la limite lorsque $t\to\infty$ des op\'erations $a\oplus_t b=\log_t(t^a+t^b)$ et $a\otimes_t b=a+b$, induites sur $\R$ en demandant que $\log_t:(\R_+^*,+,\cdot)\to (\R,\oplus_t,\otimes_t)$ soit un isomorphisme. Cette d\'eformation de structure alg\'ebrique est \'etroitement li\'ee \`a la d\'eformation de structure complexe $J_tv=\frac 1 {\log(t)}iv$, $v\in TS^1$ sur $\C^*=T^*S^1$. Nous renvoyons aux excellents textes~\cite{I,IMS} pour les bases de la g\'eom\'etrie tropicale. 

La d\'eformation de structure complexe pr\'ec\'edente a permis \`a Mikhalkin de d\'emontrer un th\'eor\`eme de correspondance entre courbes alg\'ebriques et courbes tropicales rigidifi\'ees par un nombre adapt\'e de conditions d'incidence ponctuelles~\cite{Mikhalkin}. Cette approche s'adapte au cadre r\'eel et permet de d\'ecrire en termes combinatoires les invariants de Welschinger en dimension $2$ (voir aussi~\cite{Shustin,GMS11}, ainsi que~\cite{BM} pour le cas de la dimension $3$). Comme application, mentionnons l'\'equivalence logarithmique~\cite{IKS03,IKS05} de l'invariant de Welschinger $\chi^d_{3d-1}$ et de l'invariant de Gromov-Witten $N_d$ dans le cas $X=\P^2$. Un r\'esultat similaire d'\'equivalence logarithmique est valable pour $\P^3$~: alors que $\chi^d_{2d}$ est nul si $d$ est pair, $\chi^d_{2d}$ est \'equivalent en \'echelle logarithmique \`a $N_d$ lorsque $d$ est impair~\cite{BM}. Soulignons le fait que ces r\'esultats sont obtenus pour la valeur maximale admise de $r$.

Dans~\cite{ABLdM,IKS09} les auteurs d\'emontrent des formules r\'ecursives tropicales pour calculer les invariants de Welschinger. Ces formules de type Caporaso-Harris \cite{CS98,GM07} font intervenir des invariants tropicaux \emph{relatifs}, auxquels on ne sait pas encore donner un sens en termes de courbes $J$-holomorphes. Notons au passage le lien \'etroit entre la preuve de la formule de Caporaso-Harris~\cite{CS98} et la proc\'edure d'\'etirement du cou d\'ecrite au~\S\ref{sec:SFT}.

%
%

\smallskip

\subsection{En guise de conclusion} 

Les r\'esultats que nous avons pr\'esent\'es indiquent que les invariants de Welschinger sont les bons analogues r\'eels des invariants de Gromov-Witten avec des conditions d'incidence ponctuelles. \`A la diff\'erence des invariants de Gromov-Witten, les invariants de Welschinger n'ont pas encore engendr\'e de \emph{th\'eorie} syst\'ematique comparable \`a celle de la cohomologie quantique ou encore \`a la sym\'etrie miroir. Le travail de Welschinger constitue l'un des d\'eclics essentiels de ces d\'eveloppements futurs.

\smallskip

\noindent {\bf Remerciements.} Je remercie Ilia Itenberg, Viatcheslav Kharlamov et Jean-Yves Welschinger pour leurs explications \'eclairantes. Je remercie ma famille pour son soutien.



\end{document}